# CLASSICAL SOLUTIONS TO REACTION–DIFFUSION SYSTEMS FOR HEDGING PROBLEMS WITH INTERACTING ITÔ AND POINT PROCESSES[1]

BY DIRK BECHERER AND MARTIN SCHWEIZER

*Imperial College London and ETH Zürich*

We use probabilistic methods to study classical solutions for systems of interacting semilinear parabolic partial differential equations. In a modeling framework for a financial market with interacting Itô and point processes, such PDEs are shown to provide a natural description for the solution of hedging and valuation problems for contingent claims with a recursive payoff structure.

**1. Introduction.** Reaction–diffusion systems are systems of semilinear parabolic partial differential equations which can interact in a possibly nonlinear way. They appear as models for phenomena from various areas of applications, ranging from ecological systems and biological pattern formation to chemical reactions; see Smoller (1994) for references. This article is concerned with applications to hedging and valuation problems in mathematical finance. Standard existence and uniqueness results for reaction–diffusion systems may not apply here because the coefficient functions of typical parametrizations in finance often are unbounded or do not satisfy linear growth constraints. A first contribution of this article is to address that issue by proving results on classical solutions in a fairly general context. A second contribution is an integrated treatment of contingent claims in the context of a flexible Markovian framework which incorporates new features and includes more specific models studied so far. We use the developed PDE techniques to provide results on the valuation and hedging of claims with a recursive payoff structure.

Received August 2003; revised April 2004.

[1]Supported by Deutsche Forschungsgemeinschaft via Graduiertenkolleg "Stochastic Processes and Probabilistic Analysis," TU Berlin, and the European Commission 5th Framework Program "Improving Human Potential."

The article is structured as follows. Section 2 contributes existence and uniqueness results for classical solutions of reaction–diffusion systems by showing existence of a unique fixed point of a suitable Feynman–Kac representation, provided that the underlying diffusion stays within a given domain and the coefficient functions satisfy local conditions. This first yields results under a (global) Lipschitz condition on the interaction term, and these are extended to a local Lipschitz condition by exploiting an additional monotonicity assumption on the interaction. The latter arises naturally in valuation problems from mathematical finance. Section 3 introduces for the subsequent applications a stochastic model with interacting Itô and point processes and gives a construction by a change of measure. The model consists of a system of stochastic differential equations which describes the Markovian dynamics of an Itô process $S$ and a further finite-state process $\eta$ driven by the point processes. This SDE system is nonstandard in that the driving process can itself depend on the solution, similarly as in Jacod and Protter (1982).

In Section 4, this framework is used as a model for an incomplete financial market, with $S$ describing the prices of tradable assets, for example, stock indices, while the process $\eta$ represents further (not directly tradable) sources of financial risk, for example, rating and credit events or the state of an insurance contract. An important feature is that our model allows for a mutual dependence between $S$ and $\eta$, in that the drift and volatility of $S$ can depend on the finite-state process $\eta$ while the intensities for changes of $\eta$ can in turn depend on the current value of $S$. In the context of mathematical finance, this can be seen as both a fusion and a generalization of a Markov chain modulated diffusion model of Black–Scholes type, as in Di Masi, Kabanov and Runggaldier (1994), and of the Cox process model for credit risk from Lando (1998) or the conditional Markov chain model, respectively. Another contribution is that we not only study a pure pricing approach under an a priori given pricing measure, but use a combination of valuation and hedging ideas to determine both a valuation and a locally risk-minimizing hedging strategy. Following Duffie, Schroder and Skiadas (1996), we allow for claims whose payoffs can depend not just on the state variables $S$ and $\eta$ but also on the valuation process of the claim itself. This leads to a recursive valuation problem, and it turns out that reaction–diffusion systems as in Section 2 provide a natural, convenient and constructive description in terms of PDEs for the solution to valuing and hedging problems for such claims. Section 5 discusses application examples and possible extensions, mainly with a view towards hedging of credit risk. This illustrates the flexibility of our model and results and also relates them to the existing literature on this topic which gained much interest recently; see Jeanblanc and Rutkowski (2003) for an overview and more references.

**2. Classical solutions for reaction–diffusion systems.** In this section, we use stochastic methods to derive existence and uniqueness results for *classical solutions* of interacting systems of semilinear parabolic partial differential equations (PDEs). Such systems are also known as *reaction–diffusion equations* and play in subsequent sections a key role in our solutions to various valuation and hedging problems from mathematical finance. There we consider a Markovian setting where an Itô process $S$ models the prices of the tradable assets, and further nontraded factors of risk are represented by a finite-state process $\eta$. Similarly as in the Black–Scholes model, the solutions to our valuation and hedging problems can be conveniently described via PDEs. But the nontradable factors lead to an *interacting system* of PDEs; each single PDE corresponds to a possible state of $\eta$, and the interaction between the PDEs reflects the impact from the evolution of $\eta$ on the valuation.

We first derive results for PDE systems where the interaction satisfies a global Lipschitz condition, and extend these to a type of monotonic local Lipschitz interaction. The latter is relevant for the applications to valuation and hedging problems. We strive for general assumptions on the coefficient functions which are satisfied by typical financial models.

2.1. *General framework.* Fix $m \in \mathbb{N}$, a time horizon $T \in (0, \infty)$ and a domain (open connected subset) $D$ in $\mathbb{R}^d$. For each starting point $(t, x, k) \in [0, T] \times D \times \{1, \ldots, m\}$, consider the following stochastic differential equation (SDE) in $\mathbb{R}^d$:

$$
(2.1) \quad \begin{aligned} X_t^{t,x,k} &= x \in D, \\ dX_s^{t,x,k} &= \Gamma_k(s, X_s^{t,x,k})\,ds + \sum_{j=1}^r \Sigma_{k,j}(s, X_s^{t,x,k})\,dW_s^j, \qquad s \in [t, T], \end{aligned}
$$

for continuous functions $\Gamma_k : [0, T] \times D \to \mathbb{R}^d$ and $\Sigma_{k,j} : [0, T] \times D \to \mathbb{R}^d$, $j = 1, \ldots, r$, with an $\mathbb{R}^r$-valued Brownian motion $W = (W^j)_{j=1,\ldots,r}$. We write $\Gamma_k$ and each $\Sigma_{k,j}$ as a $d \times 1$ column vector and define the matrix-valued function $\Sigma_k : [0, T] \times D \to \mathbb{R}^{d \times r}$ by $\Sigma_k^{ij} := (\Sigma_{k,j})^i$. For any $k$, $\Gamma_k$ and $\Sigma_{k,j}$, $j = 1, \ldots, r$, are assumed locally Lipschitz-continuous in $x$, uniformly in $t$:

(2.2)
> For each compact subset $\mathcal{K}$ of $D$,
> there is a constant $c = c(\mathcal{K}) < \infty$ such that
> $$|G(t, x) - G(t, y)| \le c|x - y|$$
> for all $t \in [0, T]$, $x, y \in \mathcal{K}$ and $G \in \{\Gamma_k, \Sigma_{k,1}, \ldots, \Sigma_{k,r}\}$.

By Theorem V.38 in Protter (2004), condition (2.2) implies that (2.1) has a unique (strong) solution for any given tuple $(\Omega, \mathcal{F}, \mathbb{F}, P, W)$ up to a possibly finite random explosion time. We impose the additional global and probabilistic assumption that for all $(t, x, k)$, the solution $X^{t,x,k}$ does not leave $D$ before $T$, that is,

(2.3) $$P[X_s^{t,x,k} \in D \text{ for all } s \in [t, T]] = 1.$$

This includes that $X^{t,x,k}$ does not explode to infinity:

$$P\left[\sup_{s \in [t,T]} |X_s^{t,x,k}| < \infty\right] = 1.$$

By Theorem II.5.2 of Kunita (1984), (2.2) and (2.3) imply that $X^{t,x,k}$ has a version such that

(2.4) $$(t, x, s) \mapsto X_s^{t,x,k} \text{ is } P\text{-a.s. continuous.}$$

2.2. *Fixed points of the Feynman–Kac representation (generalized solutions).* For existence and uniqueness problems of nonlinear PDEs, it is common to consider generalized solutions, namely solutions of a corresponding integral equation. These in general require and possess less regularity, and additional assumptions are needed to ensure that a solution to the integral equation is also a classical solution to the PDE. See Chapter 6.1 of Pazy (1983) for an analytic version of this approach and Freidlin (1985) for a probabilistic version.

For the PDE (2.8) that we consider below, the integral form of the stochastic approach is the well-known Feynman–Kac representation. Since the PDE is nonlinear, the solution itself appears within the expectation so that we have to look for a fixed point. To make this precise, we start with continuous functions $h: D \to \mathbb{R}^m$, $g: [0, T] \times D \times \mathbb{R}^m \to \mathbb{R}^m$ and $c: [0, T] \times D \to \mathbb{R}^m$. Given sufficient integrability, one can then define an operator $F$ on functions $v$ by

(2.5) $$(Fv)^k(t, x) := E\left[h^k(X_T^{t,x,k}) e^{\int_t^T c^k(s, X_s^{t,x,k})\, ds} + \int_t^T g^k(s, X_s^{t,x,k}, v(s, X_s^{t,x,k})) e^{\int_t^s c^k(u, X_u^{t,x,k})\, du}\, ds\right]$$

with $k = 1, \ldots, m$, $(t, x) \in [0, T] \times D$. Under suitable conditions, $F$ has a unique fixed point:

PROPOSITION 2.1. *Assume (2.2) and (2.3) hold. Let $h$, $g$ and $c$ be continuous, with $h$ and $g$ bounded, and with $c$ bounded from above in all coordinates. Then $F$ defines a mapping from the Banach space $C_b := C_b([0, T] \times$* $D, \mathbb{R}^m)$ of continuous bounded functions $v : [0,T] \times D \to \mathbb{R}^m$ into itself. Assume further that $(t, x, v) \mapsto g(t, x, v)$ is Lipschitz-continuous in $v$, uniformly in $t$ and $x$, that is, there exists $L < \infty$ such that

$$|g(t, x, v_1) - g(t, x, v_2)| \leq L|v_1 - v_2| \tag{2.6}$$

for all $t \in [0, T]$, $x \in D$ and $v_1, v_2 \in \mathbb{R}^m$.

Then $F$ is a contraction on $C_b$ with respect to the norm

$$\|v\|_\beta := \sup_{(t,x) \in [0,T] \times D} e^{-\beta(T-t)} |v(t,x)| \tag{2.7}$$

for $\beta < \infty$ large enough. In particular, $F$ has then a unique fixed point $\hat{v} \in C_b$.

PROOF. By the boundedness assumptions and (2.3), (2.4), the operator $F$ maps $C_b$ into itself. The norm (2.7) is equivalent to the usual sup-norm on $C_b$, and by assumption, $c$ has values in $(-\infty, K]^m$ for some constant $K \in [0, \infty)$. Using this and (2.6), we obtain for $v, w \in C_b$ and $\beta > 0$ that

$$e^{-\beta(T-t)} |(Fv)^k(t,x) - (Fw)^k(t,x)|$$
$$= \frac{1}{e^{\beta(T-t)}} \left| E\left[ \int_t^T (g^k(s, X_s^{t,x,k}, v(s, X_s^{t,x,k})) \right. \right.$$
$$\left. \left. - g^k(s, X_s^{t,x,k}, w(s, X_s^{t,x,k}))) e^{\int_t^s c^k(u, X_u^{t,x,k})\,du}\,ds \right] \right|$$
$$\leq \frac{e^{KT}}{e^{\beta(T-t)}} E\left[ \int_t^T |g^k(s, X_s^{t,x,k}, v(s, X_s^{t,x,k})) \right.$$
$$\left. - g^k(s, X_s^{t,x,k}, w(s, X_s^{t,x,k}))| e^{-\beta(T-s)} e^{\beta(T-s)}\,ds \right]$$
$$\leq \frac{e^{KT}}{e^{\beta(T-t)}} L \|v - w\|_\beta \int_t^T e^{\beta(T-s)}\,ds$$
$$\leq \frac{L e^{KT}}{\beta} \|v - w\|_\beta$$

for all $(t, x) \in [0, T] \times D$ and $k = 1, \ldots, m$. Thus $F$ is a contraction for $\beta > Le^{KT}$. □

2.3. *Classical solutions under local regularity.* Now define operators $\mathcal{L}^k$, $k = 1, \ldots, m$, on sufficiently smooth functions $f : [0,T] \times D \to \mathbb{R}$ by

$$(\mathcal{L}^k f)(t,x) = \sum_{i=1}^d \Gamma_k^i(t,x) \frac{\partial f}{\partial x^i}(t,x) + \frac{1}{2} \sum_{i,j=1}^d a_k^{ij}(t,x) \frac{\partial^2 f}{\partial x^i \partial x^j}(t,x)$$

with

$$a_k(t,x) = (a_k^{ij}(t,x))_{i,j=1,\ldots,d} := \Sigma_k(t,x)\Sigma_k^{\mathrm{tr}}(t,x).$$

Consider the following system of semilinear PDEs with $k=1,\ldots,m$ and boundary conditions at terminal time $T$:

$$\frac{\partial}{\partial t}v^k(t,x) + \mathcal{L}^k v^k(t,x) + c^k(t,x)v^k(t,x) + g^k(t,x,v(t,x)) = 0$$

(2.8)
$$\text{for } (t,x) \in [0,T) \times D,$$

$$v^k(T,x) = h^k(x) \qquad \text{for } x \in D.$$

These $m$ PDEs are interacting via the $g$-term which may depend on all components of $v(t,x) = (v^k(t,x))_{k=1,\ldots,m}$. Our goal is to show that the fixed point $\hat{v}$ from Proposition 2.1 is the unique bounded classical solution to (2.8). To this end we apply a Feynman–Kac type result from Heath and Schweizer (2001) that relies on classical results by Friedman (1975) and requires only local assumptions on the coefficient functions. Note that the subsequent results do not require further boundary conditions for the PDE (2.8); this is due to condition (2.3). We impose the following additional conditions on the coefficients of the SDE (2.1) and the PDE (2.8):

(2.9) There exists a sequence $(D_n)_{n\in\mathbb{N}}$ of bounded domains with closure $\bar{D}_n \subseteq D$ such that $\bigcup_{n=1}^\infty D_n = D$, each $D_n$ has a $C^2$-boundary, and for each $n$ and $k=1,\ldots,m$,

(2.10) the functions $\Gamma_k$ and $a_k = \Sigma_k \Sigma_k^{\mathrm{tr}}$ are uniformly Lipschitz-continuous on $[0,T] \times \bar{D}_n$,

(2.11) $\det a_k(t,x) \neq 0$ for all $(t,x) \in [0,T] \times D$,

(2.12) $(t,x) \mapsto c(t,x)$ is uniformly Hölder-continuous on $[0,T] \times \bar{D}_n$,

(2.13) $(t,x,v) \mapsto g(t,x,v)$ is uniformly Hölder-continuous on $[0,T] \times \bar{D}_n \times \mathbb{R}^m$.

REMARK 2.2. We aim for a classical solution, not a generalized solution (as already provided by Proposition 2.1) or a solution in some other weak sense; see, for example, Pardoux (1999) for viscosity solutions of similar PDE systems. To the best of our knowledge, the subsequent results on *classical solutions* have not been available so far under the assumptions given here.

Apart from the global condition (2.3), we only impose mild local conditions on the coefficient functions; we assume them neither bounded nor satisfying a global (linear) growth condition over the possibly unbounded domain $D$. This is crucial because such restrictive assumptions could exclude parametrizations which are typical in financial models; the only global condition (2.3) is probabilistic and must be verified on a case by case basis, for instance by means of Feller's test for explosion. An example is given in Heath and Schweizer (2001).

Let us denote by $C^{1,2}_{(b)} := C^{1,2}_{(b)}([0,T] \times D, \mathbb{R}^m)$ the spaces of continuous (bounded) functions $v: [0,T] \times D \to \mathbb{R}^m$ which are of class $C^{1,2}$ with respect to $(t,x) \in [0,T) \times D$. Note that the $C^{1,2}$-condition is imposed only on $[0,T) \times D$ while continuity is required on all of $[0,T] \times D$.

PROPOSITION 2.3. *Assume that* (2.9)–(2.13) *hold in addition to all the assumptions for Proposition* 2.1. *Then the system* (2.8) *of interacting semilinear PDEs has a unique classical solution in* $C^{1,2}_b$, *which is given by the fixed point $\hat{v}$ from Proposition* 2.1.

PROOF. Recalling that $F$ is a contraction on $C_b$ by Proposition 2.1, we first show that $F$ maps bounded functions $w$ that are locally Hölder-continuous in $(t,x)$ on $[0,T) \times D$ into $C^{1,2}_b$, and that for each such $w$, the function $v := Fw$ satisfies the following system of $m$ PDEs with terminal conditions:

$$\frac{\partial}{\partial t} v^k(t,x) + \mathcal{L}^k v^k(t,x) + c^k(t,x) v^k(t,x) + g^k(t,x,w(t,x)) = 0,$$
(2.14)
$$(t,x) \in [0,T) \times D,$$
$$v^k(T,x) = h^k(x), \qquad x \in D.$$

It is evident from the definitions of $v$ and $F$ that $v$ satisfies the terminal condition and is bounded, with a bound on $v$ that depends only on the bounds for $h$, $g$ and $c$. To prove the above two assertions, it suffices to show for any $\varepsilon > 0$ that $v$ is in $C^{1,2}([0,T-\varepsilon] \times D, \mathbb{R}^m)$ and satisfies (2.14) on $[0,T-\varepsilon] \times D$ instead of $[0,T) \times D$. So fix arbitrary $\varepsilon \in (0,T)$ and $k \in \{1,\ldots,m\}$, and let $T' := T - \varepsilon$. For any $(t,x) \in [0,T'] \times D$, conditioning on $\mathcal{F}_{T'}$ gives

$$\begin{aligned}
v^k(t,x) &= (Fw)^k(t,x) \\
&= E\bigg[E\bigg[h^k(X^{t,x,k}_T) e^{\int_t^T c^k(s, X^{t,x,k}_s)\,ds} \\
&\qquad + \int_t^T g^k(s, X^{t,x,k}_s, w(s, X^{t,x,k}_s)) e^{\int_t^s c^k(u, X^{t,x,k}_u)\,du}\,ds \bigg| \mathcal{F}_{T'} \bigg]\bigg] \\
&= E\bigg[v^k(T', X^{t,x,k}_{T'}) e^{\int_t^{T'} c^k(s, X^{t,x,k}_s)\,ds} \\
&\qquad + \int_t^{T'} g^k(s, X^{t,x,k}_s, w(s, X^{t,x,k}_s)) e^{\int_t^s c^k(u, X^{t,x,k}_u)\,du}\,ds\bigg]
\end{aligned}$$

by using the Markov property of $X^{t,x,k}$ for the last equality; compare the proof of Theorem 1 in Heath and Schweizer (2001). Using that Theorem 1 and the above representation for $v^k(t,x)$ on $[0,T') \times D$, we obtain that $(v^k)_{k=1,\ldots,m}$ is in $C^{1,2}([0,T') \times D, \mathbb{R}^m)$ and satisfies the PDE

$$\frac{\partial}{\partial t} v^k(t,x) + \mathcal{L}^k v^k(t,x) + c^k(t,x) v^k(t,x) + g^k(t,x,w(t,x)) = 0,$$
(2.15)
$$(t,x) \in [0,T') \times D,$$

if we can verify the assumptions [A1]–[A3] from Heath and Schweizer (2001). [A1] and [A2] are precisely (2.2) and (2.3). We check the list [A3'] of conditions in Heath and Schweizer (2001) whose combination implies [A3]. Conditions [A3'], [A3a'] and [A3c'] in Heath and Schweizer (2001) are exactly (2.9), (2.10) and (2.12). By Lemma 3 in Heath and Schweizer (2001), the continuity of $\Sigma_k$ in combination with (2.11) implies their condition [A3b'] ($a$ is uniform elliptic on $[0,T] \times D_n$), and $v \in C_b$ implies [A3e'] ($v$ is finite and continuous). To verify [A3d'] [i.e., (2.16) below], note that $w$ is by assumption uniformly Hölder-continuous on the compact subsets $[0,T'] \times \bar{D}_n$. Hence (2.13) implies that the composition

(2.16) $(t,x) \mapsto g(t,x,w(t,x))$ is uniformly Hölder-continuous on $[0,T'] \times \bar{D}_n$

which is [A3d']. Since $\varepsilon > 0$ was arbitrary, we conclude from (2.15) that $v = (v^k)_{k=1,\ldots,m} = Fw$ is in $C_b^{1,2}$ and satisfies the PDE (2.14) on $[0,T) \times D$.

It follows that the fixed point $\hat{v} \in C_b$ from Proposition 2.1 is approximated in $C_b$ (i.e., in the sup-norm) by a sequence $(v_n)_{n \in \mathbb{N}_0} = (Fv_{n-1})_{n \in \mathbb{N}}$ from $C_b^{1,2}$ if we choose $v_0 = w$ locally Hölder-continuous in $(t,x)$. To prove that $\hat{v}$ is an element of $C_b^{1,2} \subset C_b$ and satisfies (2.8), it suffices by the preceding argument to show that

(2.17)  $\hat{v}$ is locally Hölder-continuous in $(t,x)$ on $[0,T) \times D$.

To establish (2.17), we employ an a priori Hölder estimate for the sequence $(v_n)_{n \in \mathbb{N}}$ which is local in $(t,x)$ but uniform in $n$. Let $Q$ denote a bounded domain with $\bar{Q} \subseteq [0,T) \times D$, and let $Q'$ be some subdomain of $Q$ having a strictly positive distance to $\partial Q \cap (0,T)$, where $\partial Q$ denotes the boundary of $Q$. Then there is a $Q'$-local Hölder estimate which holds uniformly for all functions of the sequence $(v_n)$. More precisely, there exist by Theorems 6 and 7 in Chapter 4.2 of Krylov (1987) some constants $\beta \in (0,1)$ and $N < \infty$ such that the estimate

$$(2.18) \quad |v_n^k(t,x) - v_n^k(t',x')|$$
$$\leq N(\|v_n^k\|_{L^\infty(Q)} + \|g^k(\cdot,v_{n-1}(\cdot))\|_{L^{d+1}(Q)})(|x-x'| + |t-t'|^{1/2})^\beta$$

holds for all $(t,x),(t',x') \in Q'$, $n \in \mathbb{N}$ and $k = 1,\ldots,m$. This uses that (2.14) holds for $v = v_n = Fv_{n-1}$ and $w = v_{n-1}$; note also that $C_b^{1,2}(\bar{Q})$ is contained in the Sobolev space $W_{d+1}^{1,2}(Q)$ of Krylov (1987) because we can choose $Q$ to have smooth boundary and use approximations by $C^\infty$-functions. By the boundedness of $Q$ and $g$, the $L^\infty$ and $L^{d+1}$-norms of the functions $v_n^k(t,x)$ and $g^k(t,x,v_{n-1}(t,x))$, respectively, with respect to Lebesgue measure are bounded uniformly in $n \in \mathbb{N}$. In fact, we have (as already noted) a uniform bound on all $\|v_n^k\|_{L^\infty(Q)}$ in terms of $h$, $g$ and $c$, and

$$\|g^k(\cdot,v_{n-1}(\cdot))\|_{L^{d+1}(Q)} \leq \|1\|_{L^{d+1}(Q)} \|g^k\|_{L^\infty(Q)} < \infty \quad \text{for } k=1,\ldots,m.$$

Hence the first bracket on the right-hand side of (2.18) is bounded uniformly in $n$. For any such $Q'$, this yields a $Q'$-local Hölder estimate for the sequence $(v_n)_{n \in \mathbb{N}}$ that is uniform in $n$. It follows that the (uniform) limit $\hat{v} \in C_b$ is locally Hölder-continuous. This establishes (2.17).

Uniqueness follows by the usual Feynman–Kac argument. In fact, we apply Itô's formula to the process $v^k(s, X_s^{t,x,k}) \exp(\int_t^s c^k(u, X_u^{t,x,k}) \, du)$, $s \in [t,T]$, and use the PDE to show that any solution $v \in C_b^{1,2}$ to (2.8) is given by the Feynman–Kac representation and therefore a fixed point of $F$. Since the fixed point is unique, this yields $v = \hat{v}$. □

2.4. *Classical solutions under monotonic interaction.* For our applications in the context of hedging and pricing in finance, we need a solution to the PDE system (2.8) with a function $g \in C^1([0,T] \times D \times \mathbb{R}^m, \mathbb{R}^m)$ which is usually unbounded on $[0,T] \times D \times \mathbb{R}^m$, but locally bounded in $v$ in the sense that

(2.19) $g(t,x,v)$ is bounded on $[0,T] \times D \times \mathcal{K}$

for any compact subset $\mathcal{K}$ of $\mathbb{R}^m$,

and satisfies a monotonicity assumption of the following type:

There exist $K_1, K_2 \in [0,\infty)$ such that for all $k=1,\ldots,m$ we have

$$g^k(t,x,v) \leq +K_1 + K_2|v|$$

(2.20) for all $t \in [0,T]$, $x \in D$, $v \in \{w \in \mathbb{R}^m | w^k \geq w^j, \, \forall j \neq k\}$,

$$g^k(t,x,v) \geq -K_1 - K_2|v|$$

for all $t \in [0,T]$, $x \in D$, $v \in \{w \in \mathbb{R}^m | w^k \leq w^j, \, \forall j \neq k\}$.

Thus we have an upper (or lower) bound on coordinate $k$ of $g(t,x,v)$, at most linear in $|v|$, if the argument $v$ has its largest (or smallest) coordinate for index $k$. We further suppose that

(2.21) $(t,x,v) \mapsto g(t,x,v)$ is locally Lipschitz-continuous in $v$,

uniformly in $(t,x)$.

Typical examples from mathematical finance which satisfy (2.20) and (2.21) are

$$g^k(t,x,v) = \delta^k(t,x,v) + \sum_{\substack{j=1 \\ j \neq k}}^{m} \lambda^{kj}(t,x)(v^j - v^k + f^{kj}(t,x,v)),$$

(2.22)
$$k = 1, \ldots, m,$$

as in Section 4, or [as in Becherer (2004)]

$$g^k(t,x,v) = \delta^k(t,x,v) + \sum_{\substack{j=1 \\ j \neq k}}^{m} \lambda^{kj}(t,x)\frac{1}{\alpha}(e^{\alpha(v^j - v^k + f^{kj}(t,x,v))} - 1),$$

(2.23)
$$k = 1, \ldots, m,$$

with functions $\lambda^{kj} \in C_b^1([0,T] \times D, [0,\infty))$ and $\delta^k, f^{kj} \in C_b^1([0,T] \times D \times \mathbb{R}^m, \mathbb{R})$ for $k, j = 1, \ldots, m$, which are locally Lipschitz-continuous in $v$, uniformly in $(t,x)$, and with $\alpha > 0$. In Section 4, we shall work with (2.22) and replace boundedness of $\delta^k$ and $f^{kj}$ by a linear growth condition in $v$; this is still covered by (2.20).

Under the above assumptions on $g$, we cannot apply Proposition 2.3 directly since $g$ is not bounded, is not (globally) Lipschitz-continuous in $v$ and does not satisfy (2.13) in general. But we can still obtain the following result:

THEOREM 2.4. *Suppose $h$ and $c$ are continuous functions with $h$ bounded and with $c$ bounded from above in all coordinates. Assume that (2.2), (2.3) and (2.9)–(2.12) hold. Suppose $g$ is in $C^1$ and satisfies the local boundedness condition (2.19), the monotonicity condition (2.20) and the local Lipschitz condition (2.21). Then the system (2.8) of PDEs has a unique classical solution $\hat{v} \in C_b^{1,2}([0,T) \times D, \mathbb{R}^m)$, and $\hat{v}$ satisfies the Feynman–Kac representation*

$$\hat{v}^k(t,x) = E\bigg[h^k(X_T^{t,x,k})e^{\int_t^T c^k(s,X_s^{t,x,k})\,ds}$$

(2.24)
$$+ \int_t^T g^k(s, X_s^{t,x,k}, \hat{v}(s, X_s^{t,x,k}))e^{\int_t^s c^k(u,X_u^{t,x,k})\,du}\,ds\bigg]$$

*for $k = 1, \ldots, m$ and $(t,x) \in [0,T] \times D$.*

PROOF. It suffices to prove the assertion for the case where $c$ is nonpositive in all coordinates because the general case can be reduced to this by passing to the transformed function $(t,x,k) \mapsto \exp(-K(T-t))v(t,x,k)$ for a suitable constant $K \in [0,\infty)$.

Since $h$ is bounded, there is a constant $K_3 \in [0, \infty)$ such that $|h^k(x)| \leq K_3$ for all $x$ and $k$. For this constant and $K_{1,2}$ from (2.20) with $K_2$ taken with respect to the max-norm on $\mathbb{R}^m$, we first define a truncation-boundary function $\kappa : [0, T] \to [0, \infty)$ by

$$t \mapsto \kappa(t) := \begin{cases} K_3 + K_1(T-t), & \text{when } K_2 = 0, \\ K_3 e^{K_2(T-t)} + \dfrac{K_1}{K_2}(e^{K_2(T-t)} - 1), & \text{when } K_2 > 0, \end{cases}$$

and then $\tilde{g} : [0, T] \times D \times \mathbb{R}^m \to \mathbb{R}^m$ by truncating the third argument:

$$\tilde{g}(t, x, v) := g(t, x, (\min(\max(v^k, -\kappa(t)), +\kappa(t)))_{k=1,\ldots,m}).$$

Then $\tilde{g}$ is bounded by (2.19), and Lipschitz-continuous on $[0, T] \times \bar{D}_n \times \mathbb{R}^m$ for every $n$ since $g$ is in $C^1$. By (2.21) the function $\tilde{g}$ is also Lipschitz in $v$, uniformly for $(t, x) \in [0, T] \times D$, and so Proposition 2.3 yields a unique bounded solution $\hat{v}$ for the PDE (2.8) with $\tilde{g}$ instead of $g$. Moreover, $\hat{v}$ is the fixed point of $F$ defined with $\tilde{g}$ instead of $g$ in (2.5). We show below that

(2.25) $\quad |\hat{v}^k(t, x)| \leq \kappa(t) \qquad$ for $(t, x) \in [0, T] \times D$ and $k = 1, \ldots, m$.

Admitting this result for the moment, we get $\tilde{g}(t, x, \hat{v}(t, x)) = g(t, x, \hat{v}(t, x))$ for all $(t, x) \in [0, T] \times D$ by the definition of $\tilde{g}$. Hence, $\hat{v}$ also solves the PDE (2.8) with $g$ instead of $\tilde{g}$ and satisfies (2.24). To see that $\hat{v}$ is the unique bounded solution to (2.8), let $w$ denote another bounded solution. By taking $K_3$ larger if necessary, we can assume that $|w^k(t, x)| \leq K_3 \leq \kappa(t)$ for all $k, t, x$. Then both $\hat{v}$ and $w$ solve (2.8) not only with $g$ but also with $\tilde{g}$, and this implies $\hat{v} = w$ by the uniqueness in Proposition 2.3 applied for $\tilde{g}$.

To finish the proof, it remains to establish (2.25). Fix arbitrary $(t, x, k) \in [0, T] \times D \times \{1, \ldots, m\}$ and define the stopping time

$$\tau := \inf\{s \in [t, T] | \hat{v}^k(s, X_s^{t,x,k}) < \kappa(s)\} \wedge T.$$

Then

$$\hat{v}^k(s, X_s^{t,x,k}(\omega)) \geq \kappa(s) \qquad \text{for all } (\omega, s) \in [\![t, \tau[\![$$

while $\hat{v}^k(\tau, X_\tau^{t,x,k}) \leq \kappa(\tau)$; in fact, we have equality for $\tau < T$ and inequality for $\tau = T$ since $\hat{v}^k(T, \cdot) = h^k(\cdot) \leq K_3$. Hence, the definition of $\tilde{g}$ and property (2.20) of $g$ imply

(2.26) $\tilde{g}^k(s, X_s^{t,x,k}(\omega), \hat{v}(s, X_s^{t,x,k}(\omega))) \leq K_1 + K_2\kappa(s) \qquad$ for $(\omega, s) \in [\![t, \tau[\![$

and therefore by using $c \leq 0$ and integrating $\kappa$

$$\hat{v}^k(t, x) = E\bigg[E\bigg[h^k(X_T^{t,x,k})e^{\int_t^T c^k(s, X_s^{t,x,k})\,ds} \\ + \int_t^T \tilde{g}^k(s, X_s^{t,x,k}, \hat{v}(s, X_s^{t,x,k}))e^{\int_t^s c^k(u, X_u^{t,x,k})\,du}\,ds\bigg|\mathcal{F}_\tau\bigg]\bigg]$$

$$= E\bigg[\hat{v}^k(\tau, X_\tau^{t,x,k})e^{\int_t^\tau c^k(s, X_s^{t,x,k})\,ds}$$

$$+ \int_t^\tau \tilde{g}^k(s, X_s^{t,x,k}, \hat{v}(s, X_s^{t,x,k}))e^{\int_t^s c^k(u, X_u^{t,x,k})\,du}\,ds\bigg]$$

$$\leq E\bigg[\kappa(\tau) + \int_t^\tau (K_1 + K_2\kappa(s))\,ds\bigg]$$

$$= E[\kappa(\tau) + (\kappa(t) - \kappa(\tau))]$$

$$= \kappa(t),$$

where the second equality uses the strong Markov property of $X^{t,x,k}$; see Heath and Schweizer (2001). This gives the upper bound in (2.25), and the lower bound is proved in the same way by using $\tau := \inf\{s \in [t, T] | \hat{v}^k(s, X_s^{t,x,k}) > -\kappa(s)\} \wedge T$. This completes the proof. □

**3. A model with interacting Itô and point processes.** In this section, we introduce a flexible Markovian model for an incomplete financial market, give a rigorous construction for it and provide some further properties. This is used in the next section to derive explicit and constructive results for various hedging and valuation approaches in terms of reaction–diffusion systems. The entire model is given by a system of stochastic differential equations (SDEs) of two types. The price process $S$ of the tradable risky assets is modeled by an Itô process. In addition, there are nontradable factors of uncertainty and risk which are represented by a finite-state process $\eta$ driven by a point process. A distinctive feature of our model is that it permits mutual dependences between $S$ and $\eta$. More precisely, both $S$ and $\eta$ enter the coefficients of the SDE for the dynamics of $S$, and at the same time, the intensities controlling the jumps of $\eta$ depend on the current value of $S$.

3.1. *Model setup and assumptions.* All modeling in the sequel takes place on some filtered probability space $(\Omega, \mathcal{F}, \mathbb{F}, P)$ with a filtration $\mathbb{F} = (\mathcal{F}_t)_{t \in [0,T]}$ satisfying the usual conditions and a trivial $\sigma$-field $\mathcal{F}_0$. All semimartingales are taken to have right-continuous paths with left limits.

We start with $m \in \mathbb{N}$ and a domain $D$ in $\mathbb{R}^d$ satisfying (2.9), for example, $D = \mathbb{R}^d$ or $D = (0, \infty)^d$. Let $(S, \eta)$ be a solution of the following system of SDEs with values in $D \times \{1, \ldots, m\}$:

(3.1) $\quad S_0 \in D, \qquad dS_t = \Gamma(t, S_t, \eta_{t-})\,dt + \Sigma(t, S_t, \eta_{t-})\,dW_t,$

(3.2) $\quad \eta_0 \in \{1, \ldots, m\}, \qquad d\eta_t = \sum_{k,j=1}^m (j-k)I_{\{k\}}(\eta_{t-})\,dN_t^{kj},$

where $\Gamma : [0,T] \times D \times \{1,\ldots,m\} \to \mathbb{R}^d$ and $\Sigma : [0,T] \times D \times \{1,\ldots,m\} \to \mathbb{R}^{d \times r}$ are $C^1$ with respect to $(t,x) \in [0,T] \times D$, $W = (W^i)_{i=1,\ldots,r}$ is an $\mathbb{R}^r$-valued $(P,\mathbb{F})$-Brownian motion and $N = (N^{kj})_{k,j=1,\ldots,m}$ is a multivariate $\mathbb{F}$-adapted point process such that

(3.3) $\quad (N_t^{kj})$ has $(P,\mathbb{F})$-intensity $\lambda^{kj}(t,S_t) \quad$ for $k,j = 1,\ldots,m$

with bounded $C^1$ functions $\lambda^{kj} : [0,T] \times D \to [0,\infty)$. Note that the process counting the jumps of $\eta$ from state $k$ to $j$ is not $N^{kj}$, but $\int I_{\{k\}}(\eta_-)\, dN^{kj}$. If $D \subseteq (0,\infty)^d$, one can rewrite (3.1) as a generalized Black–Scholes model; denoting $dS/S = (dS^i/S^i)_{i=1,\ldots,d}$, $\gamma(t,x,k) = \mathrm{diag}(1/x^i)_{i=1,\ldots,d}\Gamma(t,x,k)$ and $\sigma(t,x,k) = \mathrm{diag}(1/x^i)_{i=1,\ldots,d}\Sigma(t,x,k)$, we have

(3.4) $$\frac{dS_t}{S_t} = \gamma(t,S_t,\eta_{t-})\,dt + \sigma(t,S_t,\eta_{t-})\,dW_t.$$

The model (3.1)–(3.3) is a nonstandard SDE system because of its dependence structure. The coefficients in the SDE (3.1) for $S$ contain $\eta$, and the intensities in (3.3) of the point process $N$ driving $\eta$ depend in turn on $S$. We shall comment below on construction and properties.

To apply the PDE results from Section 2, we need further regularity assumptions on the coefficients of the SDE (3.1). Because $(t,x) \mapsto G(t,x,k)$ is $C^1$ on $[0,T] \times D$ for $G \in \{\Gamma,\Sigma\}$ and any $k$, the mappings $x \mapsto G(t,x,k)$ are locally Lipschitz-continuous in $x$, uniformly in $t$ and $k$. This implies as in Section 2 that there is a unique strong solution $X^{t,x,k}$ to the SDE

(3.5) $\quad \begin{aligned} X_t^{t,x,k} &= x \in D, \\ dX_s^{t,x,k} &= \Gamma(s,X_s^{t,x,k},k)\,ds + \Sigma(s,X_s^{t,x,k},k)\,dW_s, \qquad s \in [t,T], \end{aligned}$

for any $(t,x,k) \in [0,T] \times D \times \{1,\ldots,m\}$ up to a possibly finite random explosion time. As in Section 2.1 [for $\Gamma_k(t,x) := \Gamma(t,x,k)$ and $\Sigma_k(t,x) := \Sigma(t,x,k)$ there], we assume again that $X^{t,x,k}$ does not leave $D$ during $[0,T]$; that is, we suppose

(3.6) $\quad P[X_s^{t,x,k} \in D \text{ for all } s \in [t,T]] = 1 \qquad$ for any $t,x,k$.

REMARK 3.1. Intuitively, the SDEs (3.5) are related to (3.1) in the sense that $S$ could be constructed successively from one jump time of $\eta$ to the next by "pasting together" appropriate solutions to (3.5), using techniques similar to those known from the construction of finite-state Markov processes or Cox processes. But such a construction of $(S,\eta)$ becomes tedious in details, and the construction via a change of measure described in Section 3.3 appears more convenient.

3.2. *Markov property and uniqueness in distribution.* A standard way to show the Markov property is to prove uniqueness of a corresponding (time-inhomogeneous) martingale problem. We give here a direct argument which is similar in spirit. For a horizon $T' \in [0, T]$ and a function $h \in C_b(D \times \{1, \ldots, m\}, \mathbb{R})$, we consider the PDE system

$$
\begin{aligned}
0 = v_t(t,x,k) + \Gamma(t,x,k)\operatorname{grad}_x v(t,x,k) + \tfrac{1}{2} \sum_{i,j=1}^d a^{ij}(t,x,k) v_{x^i x^j}(t,x,k) \\
+ \sum_{\substack{j=1 \\ j \neq k}}^m \lambda^{kj}(t,x)(v(t,x,j) - v(t,x,k)), \qquad (t,x) \in [0,T') \times D,
\end{aligned}
\tag{3.7}
$$

for $k \in \{1, \ldots, m\}$, with $(a^{ij})_{i,j=1,\ldots,d} = a := \Sigma \Sigma^{\mathrm{tr}}$, and terminal conditions

$$v(T',x,k) = h(x,k), \qquad x \in D. \tag{3.8}$$

For brevity, we use subscripts for the partial derivatives of $v$. By Theorem 2.4 there is a unique bounded classical solution $v \in C_b^{1,2}([0,T') \times D \times \{1,\ldots,m\}, \mathbb{R})$ to (3.7) and (3.8) for any given $T'$ and $h$. The essential martingale argument for Proposition 3.3 is

LEMMA 3.2. *For $v$ given as above, the process $v(t, S_t, \eta_t)$, $t \in [0, T']$, is a martingale.*

The proof is mainly an application of Itô's formula and given in the Appendix. In the same way, another application of Itô's formula yields that $(S, \eta)$ solves the following martingale problem: For any continuous function $f(x,k)$ on $D \times \{1, \ldots, m\}$ with compact support that is of class $C^2$ in $x$, the process

$$f(S_t, \eta_t) - f(S_0, \eta_0) - \int_0^t \mathcal{A}_s f(S_s, \eta_s)\,ds, \qquad t \in [0,T],$$

is a martingale, with the operators $\mathcal{A}_s$ being given by

$$
\begin{aligned}
\mathcal{A}_s f(x,k) = \Gamma(s,x,k)\operatorname{grad}_x f(x,k) + \tfrac{1}{2} \sum_{i,j=1}^d a^{ij}(s,x,k) f_{x^i x^j}(x,k) \\
+ \sum_{\substack{j=1 \\ j \neq k}}^m \lambda^{kj}(s,x)(f(x,j) - f(x,k)).
\end{aligned}
\tag{3.9}
$$

PROPOSITION 3.3. *$(S_t, \eta_t)$, $t \in [0,T]$, is a (time-inhomogeneous) Markov process with respect to $P$ and $\mathbb{F}$. Its distribution is uniquely determined by the SDE system* (3.1)–(3.3).

PROOF. For any $h \in C_b(D \times \{1, \ldots, m\}, \mathbb{R})$ and $T' \in [0, T]$ there is a unique bounded classical solution $v$ to the PDE (3.7) with terminal condition (3.8). By Lemma 3.2,

$$E[h(S_{T'}, \eta_{T'})|\mathcal{F}_t] = E[v(T', S_{T'}, \eta_{T'})|\mathcal{F}_t] = v(t, S_t, \eta_t) \qquad \text{for } 0 \leq t \leq T'.$$

This establishes the Markov property of $(S, \eta)$. In particular, $t = 0$ gives $E[h(S_{T'}, \eta_{T'})] = v(S_0, \eta_0)$ and so the one-dimensional marginal distributions of the process are unique. To show uniqueness of the finite-dimensional distributions by induction, let $h^1, \ldots, h^{n+1}$ be arbitrary continuous bounded functions. For any times $t_1 \leq \cdots \leq t_{n+1} \leq T$, conditioning on $\mathcal{F}_{t_n}$ gives

$$(3.10) \qquad E\left[\prod_{i=1}^{n+1} h^i(S_{t_i}, \eta_{t_i})\right] = E\left[\left(\prod_{i=1}^{n} h^i(S_{t_i}, \eta_{t_i})\right) v(S_{t_n}, \eta_{t_n})\right],$$

where $v$ denotes the solution to the PDE (3.7) and (3.8) with $h := h^{n+1}$ and $T' := t_{n+1}$. Since the right-hand side of (3.10) is determined by the $n$-dimensional distributions, the claim follows. $\square$

3.3. *Construction by a change of measure.* At first sight, the mutual dependences in (3.1)–(3.3) might seem to make such models difficult to construct since we face a nonstandard SDE system where the solution $(S, \eta)$ also affects one part $N$ of the driving process. But the problem can be reduced to the special case where $N = (N^{kj})$ is a standard multivariate point process; then $\eta$ becomes an autonomous process and $S$ is well defined by (3.1). From here, the desired $(t, S)$-dependent intensities for $N$ can then be constructed by a suitable change of measure.

More precisely, we start with a filtered probability space $(\Omega, \mathcal{F}', \mathbb{F}', P')$ carrying an $r$-dimensional $(P', \mathbb{F}')$-Brownian motion $W = (W^i)_{i=1,\ldots,r}$ and a multivariate $\mathbb{F}'$-adapted point process $N = (N^{kj})_{k,j=1,\ldots,m}$ with constant $(P', \mathbb{F}')$-intensity 1 for any $k$ and $j$. In other words,

(3.11) $\qquad N^{kj}, k, j = 1, \ldots, m$, are independent standard Poisson processes under $P'$.

We assume that $\mathcal{F}'_0$ is trivial, $\mathcal{F}'_T = \mathcal{F}'$ and $\mathbb{F}'$ satisfies the usual conditions. Then (3.2) defines an autonomous process $\eta$. Given this process, there is a solution $S$ to (3.1) under suitable assumptions on the coefficients; simple examples are given in the following

EXAMPLE 3.4. Provided (3.11) holds, there exists a solution $(S, \eta)$ to (3.1), (3.2), and the solutions $X^{t,x,k}$ to (3.5) satisfy condition (3.6) in the following cases [cf. Becherer (2004)]:

(i) $D = (0, \infty)^d$, and $\gamma, \sigma$ in (3.4) are continuous functions, depending only on $(t, k)$ but not on $x$. In this case $X^{t,x,k}$ and $S$ can even be written explicitly as stochastic exponentials.

(ii) $D = \mathbb{R}^d$, and $\Gamma, \Sigma$ are Lipschitz-continuous in $x$, uniformly in $t$, for any $k$.

For the general case, we define a probability measure $P \ll P'$ by

$$(3.12) \qquad dP := \mathcal{E}\left(\sum_{k,j=1,\ldots,m} \int (\lambda^{kj}(t, S_t) - 1)(dN_t^{kj} - dt)\right)_T dP',$$

where the density is chosen to give $N$ the $(P, \mathbb{F}')$-intensities (3.3); see Chapter VI.2 in Brémaud (1981). By Girsanov's theorem, $W$ is a local $(P, \mathbb{F}')$-martingale whose covariance process $\langle W \rangle$ is the same under $P'$ and $P$ since it can be computed pathwise, and therefore $W$ is also a $(P, \mathbb{F}')$-Brownian motion. Finally, if $(\Omega, \mathcal{F}, \mathbb{F}, P)$ is the standard $P$-completion of $(\Omega, \mathcal{F}', \mathbb{F}', P)$, one can check that $\mathbb{F}$ satisfies the usual conditions under $P$. With respect to $(P, \mathbb{F})$ we then have that $W$ is a Brownian motion, $N$ is a multivariate point process with the desired intensities (3.3) and $(S, \eta)$ solves (3.1) and (3.2).

This change-of-measure construction extends an argument from Kusuoka (1999) on how to construct two point processes with mutually dependent intensities. Alternatively, one could infer existence (but not uniqueness) of a solution to the SDE system (3.1)–(3.3) from results by Jacod and Protter (1982); they constructed a solution for an SDE where the characteristics of the driving process depend on the solution process by transforming the problem to an SDE of ordinary type on a suitably enlarged probability space.

**4. Hedging and valuation of integrated risks.** This section presents an approach for valuing and hedging a general class of contingent claims with recursive payoff structure in the model introduced in the last section, and shows how option values and hedging strategies for this framework can be obtained from solutions of suitable PDE systems of the type studied in Section 2. One important feature of the claims we consider is that they can be specified implicitly, in the sense that their own value may influence the payoffs they deliver. A typical example is the pricing of a defaultable zero coupon bond with fractional recovery when the recovered amount depends on the pre-default value of the bond. Such claims lead to a fixed-point problem since their value depends on their payoff structure, which in turn depends on their value. We show how this fits naturally into the framework developed so far and leads, in comparison to the general setting of Duffie, Schroder and Skiadas (1996), to more explicit results in terms of PDEs in our setting. In addition to a pure pricing approach, we also offer a combination of valuing and hedging ideas.

4.1. *Formulation of the problem.* Our model of the financial market is given by the process $(S, \eta)$ from (3.1)–(3.3) with the assumptions from Section 3.1. We think of $S$ as the prices of $d$ risky assets (discounted, i.e., expressed in some tradable numeraire) and of $\eta$ as representing some nontradable risk factors in the market. Recall that $\Gamma$, $\Sigma$ and $\lambda^{kj}$ are $C^1$ in $(t, x)$ and all $\lambda^{kj}$ are bounded. In addition, we suppose that the market price of risk function

$$(4.1)\quad \Phi := \Sigma^{\mathrm{tr}}(\Sigma\Sigma^{\mathrm{tr}})^{-1}\Gamma \text{ exists and is bounded on } [0, T] \times D \times \{1, \ldots, m\}.$$

This implies that $\widehat{Z} = \mathcal{E}(-\int \Phi(t, S_t, \eta_{t-})\, dW_t)$ is in the Hardy martingale space $\mathcal{H}^p(P)$ for any $p \in [1, \infty)$ and so $d\widehat{P} = \widehat{Z}_T\, dP$ defines an equivalent local martingale measure (ELMM) for $S$, the so-called minimal ELMM. By Girsanov's theorem, the dynamics of $S$ under $\widehat{P}$ are given by

$$S_0 \in D, \qquad dS_t = \Sigma(t, S_t, \eta_{t-})\, d\widehat{W}_t$$

for a $\widehat{P}$-Brownian motion $\widehat{W} = W + \int \Phi(t, S_t, \eta_{t-})\, dt$, while the dynamics of $\eta$ and $N$ under $\widehat{P}$ are the same as under $P$ and given by (3.2) and (3.3). Hence $(S, \eta)$ is also Markov under $\widehat{P}$, and

$$(4.2)\qquad M^{kj} := N^{kj} - \int \lambda^{kj}(t, S_t)\, dt, \qquad k, j = 1, \ldots, m,$$

are martingales in $\mathcal{H}^p$ for any $p \in [1, \infty)$, under both $P$ and $\widehat{P}$; in fact, since $[M^{kj}] = N^{kj}$ and the intensities $\lambda^{kj}$ are bounded, a simple time change argument via Theorem II.16 in Brémaud (1981) shows that $[M^{kj}]_T$ has even all exponential moments.

REMARK 4.1. Except in trivial cases, there are typically many other ELMMs beside $\widehat{P}$ so that our financial market is incomplete. Note also that $\widehat{P}$ and $P$ coincide for $\Gamma = 0$, that is, when the original probability measure is already a local martingale measure for $S$. This automatically happens if one starts modeling under a pricing measure, as practitioners often do.

The financial contracts we consider are specified by functions $h: D \times \{1, \ldots, m\} \to \mathbb{R}$, $\delta: [0, T] \times D \times \{1, \ldots, m\} \times \mathbb{R} \to \mathbb{R}$ and $f^{kj}: [0, T] \times D \times \mathbb{R} \to \mathbb{R}$, $k, j = 1, \ldots, m$, with the following interpretations. [Note a slight change of notation in comparison to (2.22). Technically, we shall find that the function $v$ is given by a PDE system as in Section 2, where the interaction $g$ has the form (2.22) with $\delta^k$ and $f^{kj}$ depending on the argument $v$ only via the $k$th coordinate $v^k$; this is sufficient for later applications and allows to simplify notation in the sequel. Formally, the relation to (2.22) is given by $\tilde{\delta}^k(t, x, v) = \delta(t, x, k, v^k)$ and $\tilde{f}^{kj}(t, x, v) = f^{kj}(t, x, v^k)$ when $\tilde{\delta}^k$, $\tilde{f}^{kj}$ denote the functions from Section 2.] $h$ describes a final payoff at time $T$ of amount $h(x, k)$ if $S_T = x$ and $\eta_T = k$; $\delta$ specifies a rate for payments made continuously in time (e.g., dividends) and $f^{kj}$ describes a lump sum payment that falls due whenever the state of $\eta$ changes from $k$ to $j$. A typical example is given by a life insurance contract where $\eta$ could describe the state of health of the insured person. The total payoff up to time $T$ from a triple $(h, \delta, f^{kj})$ is

$$
\begin{aligned}
H = h(S_T, \eta_T) &+ \int_0^T \delta(t, S_t, \eta_{t-}, v(t, S_t, \eta_{t-})) \, dt \\
&+ \int_0^T \sum_{\substack{k,j=1 \\ k \neq j}}^m f^{kj}(t, S_t, v(t, S_t, \eta_{t-})) I_{\{k\}}(\eta_{t-}) \, dN_t^{kj}.
\end{aligned}
\tag{4.3}
$$

Technically, we always suppose that $h$, $\delta$ and $f^{kj}$ are continuous; $h$ is bounded; there is some constant $K < \infty$ such that $|\delta(t, x, k, v)|$ and $|f^{kj}(t, x, v)|$ are bounded by $K(1 + |v|)$ for all $t$, $x$, $k$, $j$; and each $\delta(\cdot, \cdot, k, \cdot)$ and all $f^{kj}$ are $C^1$ in $(t, x, v)$ and moreover locally Lipschitz-continuous in $v$, uniformly in $(t, x)$. By (3.2), one can rewrite the integral with respect to $N$ in (4.3) as a sum $\sum_{t \in (0,T]} \sum_{k \neq j} f^{kj}(\cdots) I_{\{k\}}(\eta_{t-}) I_{\{j\}}(\eta_t)$; this shows that the payoff involves only the random processes $S$ and $\eta$. Similar remarks apply at several points in the sequel, compare (4.4), (4.6) and (4.7), (4.13) and (4.14), where we use the integral notation for technical convenience.

A closer look at (4.3) shows that $H$ has an unusual feature. It is not (yet) well defined as a random variable because the extra argument $v(t, S_t, \eta_{t-})$ in $\delta$ and $f^{kj}$ has not yet been specified. The idea is that many intermediate payments depend not only on the current state of the underlying $S$ and the factor $\eta$, but also on the value of the contract itself. In life insurance, for instance, the payoff at death might be set equal to the current reserves; another example occurs with fractional recovery of a defaultable bond. However, the function $v$ that should give the current value of the contract still needs to be determined by some argument.

A first possible idea for finding $v$ is a *pure pricing* approach. We fix some probability measure and axiomatically define the value process of a (future) payment stream as the conditional expectations of all future payments under that measure. This is a common procedure in, for example, credit risk valuation problems. Since the modeling takes place under the valuation measure, we may take $P$ as that measure and thus assume that $P$ is an ELMM and $\Gamma \equiv 0$. By the Markov nature of the model, the valuation should be determined via a valuation function $v$, and in view of the recursive payoff structure we should like to write the value at time $t$ of future payoffs as

$$v(t, S_t, \eta_t) = E\left[h(S_T, \eta_T) + \int_t^T \delta(u, S_u, \eta_{u-}, v(u, S_u, \eta_{u-})) \, du\right.$$

$$(4.4) \qquad + \int_t^T \sum_{\substack{k,j=1 \\ k \neq j}}^m f^{kj}(u, S_u, v(u, S_u, \eta_{u-})) I_{\{k\}}(\eta_{u-}) \, dN_u^{kj} \bigg| \mathcal{F}_t \bigg],$$

$$t \in [0, T].$$

But since the right-hand side itself contains the function $v$, existence and uniqueness of a solution to (4.4) is a fixed-point problem; this can be viewed as a variant of the recursive valuation approach by Duffie, Schroder and Skiadas (1996).

An alternative to a pure pricing approach is a combination of *valuing and hedging* ideas. Loosely speaking, the general goal there is to find at the same time a stochastic process $\theta$ and a function $v$ such that the trading strategy $\theta$ provides a "good" hedge against the payoff $H$ in (4.3) and the value process of the strategy is related to the valuation function $v$. To make this more precise, we use here as criterion for the quality of a hedge the concept of local risk-minimization. This goes back to Schweizer (1991); the idea is to find a not necessarily self-financing strategy whose final value coincides with the sum of all payments to be hedged, and whose cumulative costs over time have in a suitable sense minimal quadratic fluctuations on average under $P$. In mathematical terms, this problem can be reformulated as finding a decomposition of the form

$$(4.5) \qquad H = H_0 + \int_0^T \theta_u \, dS_u + L_T$$

with a constant $H_0$ and a $\widehat{P}$-martingale $L$ that is $\widehat{P}$-orthogonal to $S$, where $\widehat{P}$ is the minimal ELMM introduced above. The trading strategy which holds (at any time $t$) $\theta_t^i$ units of risky asset $i = 1, \ldots, d$ and the amount $\theta_t^0 = H_0 + L_t + \int_0^t \theta_u \, dS_u - \theta_t S_t$ in the numeraire used for discounting is then locally risk-minimizing for $H$ under $P$, and one can deduce from (4.5) that the resulting valuation process $V = H_0 + \int \theta \, dS + L$ is the conditional expectation of the payoff $H$ under $\widehat{P}$. Up to integrability issues, this follows from the results in Föllmer and Schweizer (1991) since $S$ is continuous; see Theorem 3.5 in Schweizer (2001). When $P = \widehat{P}$, the above strategy is even (globally) risk-minimizing in the sense of Föllmer and Sondermann (1986). But as with (4.4), the decomposition (4.5) cannot be obtained in a standard way because $H$ involves the value function $v$ which in turn depends via $V$ on the decomposition.

REMARK 4.2. Risk-minimization is one among several hedging approaches for incomplete markets and (despite some drawbacks) typically leads to comparably constructive solutions. Since one motivation for this section comes from credit risk problems where hardly any (constructive) results for hedging complex payoffs under incompleteness seem available so far, studying risk-minimization is a natural first step to make. For results from another (nonlinear) hedging approach in a similar model but for simpler (nonrecursive) payoffs, see Becherer (2004).

Although the reasoning behind the two approaches is quite different and the pure pricing approach postulates the special case $P = \widehat{P}$, both valuations have the same mathematical structure. Both are expectations of future payments under the measure $\widehat{P}$, and we shall see that the valuation function $v$ for both is indeed determined by the recursion formula

$$
\begin{aligned}
v(t, S_t, \eta_t) = \widehat{E}\bigg[ & h(S_T, \eta_T) + \int_t^T \delta(u, S_u, \eta_{u-}, v(u, S_u, \eta_{u-}))\, du \\
& + \int_t^T \sum_{\substack{k,j=1 \\ k \neq j}}^m f^{kj}(u, S_u, v(u, S_u, \eta_{u-})) I_{\{k\}}(\eta_{u-})\, dN_u^{kj} \Big| \mathcal{F}_t \bigg],
\end{aligned}
\tag{4.6}
$$

$$t \in [0, T],$$

where $\widehat{E}$ denotes expectation under the minimal ELMM $\widehat{P}$. To prepare for this, we first establish uniqueness of the corresponding valuation process in Lemma 4.3, whose proof is relegated to the Appendix. Existence will be a by-product of the main results in the next section which describe the corresponding valuation function $v$ and the decomposition (4.5).

LEMMA 4.3. *For any payoff triple $(h, \delta, f^{kj})$ and probability measure $Q \in \{P, \widehat{P}\}$, there exists at most one bounded semimartingale $V$ such that*

$$
\begin{aligned}
V_t = E_Q\bigg[ & h(S_T, \eta_T) + \int_t^T \delta(u, S_u, \eta_{u-}, V_{u-})\, du \\
& + \int_t^T \sum_{\substack{k,j=1 \\ k \neq j}}^m f^{kj}(u, S_u, V_{u-}) I_{\{k\}}(\eta_{u-})\, dN_u^{kj} \Big| \mathcal{F}_t \bigg], \qquad t \in [0, T].
\end{aligned}
\tag{4.7}
$$

4.2. *Solution via interacting PDE systems.* For a given payoff triple $(h, \delta, f^{kj})$, consider the following system of PDEs for a function $v(t, x, k)$ on $[0, T] \times D \times \{1, \ldots, m\}$: For each $k$,

$$
\begin{aligned}
0 = {} & v_t(t, x, k) + \tfrac{1}{2} \sum_{i,j=1}^d a^{ij}(t, x, k) v_{x^i x^j}(t, x, k) + \delta(t, x, k, v(t, x, k)) \\
& + \sum_{\substack{j=1 \\ j \neq k}}^m \lambda^{kj}(t, x)(v(t, x, j) - v(t, x, k) + f^{kj}(t, x, v(t, x, k))),
\end{aligned}
\tag{4.8}
$$

$$(t, x) \in [0, T) \times D,$$

and

(4.9) $\qquad v(T, x, k) = h(x, k), \qquad x \in D,$

with $(a^{ij})_{i,j=1,\ldots,d} = a := \Sigma \Sigma^{\mathrm{tr}}$. Note that (4.8) is the special case of (2.8) where $c \equiv 0$, $\Gamma \equiv 0$ and $g$ is given by (2.22). We want to conclude from Theorem 2.4 that (4.8) and (4.9) has a unique solution in $C_b^{1,2}([0,T) \times D \times \{1, \ldots, m\}, \mathbb{R})$ for any claim $H$ of the form (4.3), and there is one point where we must take care. Assumption (3.6) from Section 3 guarantees non-explosion under $P$ for the solution of (3.5) or (2.1) with drift $\Gamma$, whereas the present application of Theorem 2.4 requires this for the solution with drift 0. However, this is also true and can be verified via Girsanov's theorem by a change to an equivalent measure under which there is no drift.

THEOREM 4.4. *Under the assumptions from Section 4.1, let* $v \in C_b^{1,2}([0,T) \times D \times \{1, \ldots, m\}, \mathbb{R})$ *be the solution of* (4.8), (4.9) *corresponding to the triple* $(h, \delta, f^{kj})$. *Then the payoff* $H$ *in* (4.3) *admits a decomposition* (4.5) *with*

$$H_0 = v(0, S_0, \eta_0),$$

(4.10) $\quad \theta_t = (\theta_t^i)_{i=1,\ldots,d} = \mathrm{grad}_x v(t, S_t, \eta_{t-}) \qquad \text{for } t \in [0, T) \quad \text{and}$

$$L_t = \int_0^t \sum_{\substack{k,j=1 \\ k \neq j}}^m (v(u, S_u, j) - v(u, S_u, k)$$

$$+ f^{kj}(u, S_u, v(u, S_u, \eta_{u-})))I_{\{k\}}(\eta_{u-}) \, dM_u^{kj},$$

*for* $t \in [0, T]$. *The process* $L$ *is a martingale in* $\mathcal{H}^p$ *under both* $P$ *and* $\widehat{P}$, *and* $\langle \int \theta \, dS \rangle_T$ *is in* $L^{p/2}(P)$ *for every* $p \in [1, \infty)$. *Moreover,* $L$ *is a BMO-martingale under both* $P$ *and* $\widehat{P}$, *and* $\int \theta \, dS$ *with* $\theta$ *from* (4.10) *is in BMO under* $\widehat{P}$.

PROOF. (i) Applying Itô's formula to $v(t, S_t, \eta_t)$ yields

$$dv(t, S_t, \eta_t) = \mathrm{grad}_x v(t, S_t, \eta_{t-}) \, dS_t$$

$$+ \left( v_t + \tfrac{1}{2} \sum_{i,j=1}^d a^{ij} v_{x^i x^j} \right)(t, S_t, \eta_{t-}) \, dt$$

$$+ \sum_{\substack{k,j=1 \\ k \neq j}}^m (v(t, S_t, j) - v(t, S_t, k))I_{\{k\}}(\eta_{t-}) \, dN_t^{kj}, \qquad t \in [0, T).$$

Substituting $dN_t^{kj} = dM_t^{kj} + \lambda^{kj}(t, S_t)\, dt$, one can then use the PDE (4.8) to obtain

$$
\begin{aligned}
dv(t, S_t, \eta_t) = {} & \theta_t\, dS_t + dL_t - \delta(t, S_t, \eta_{t-}, v(t, S_t, \eta_{t-}))\, dt \\
& - \sum_{\substack{k,j=1 \\ k \neq j}}^{m} f^{kj}(t, S_t, v(t, S_t, \eta_{t-})) I_{\{k\}}(\eta_{t-})\, dN_t^{kj}
\end{aligned}
\quad (4.11)
$$

for $t \in [0, T)$.

The next part of the proof shows in particular that $\int \theta\, dS$ is well defined on $[0, T]$ so that (4.11) extends from $[0, T)$ to all of $[0, T]$ as both sides are a.s. continuous at time $T$. Admitting this for the moment, we can integrate (4.11) from 0 to $T$ and use (4.9) to conclude by comparison with (4.3) that the ingredients in (4.10) indeed yield the decomposition (4.5).

(ii) To prove the desired integrability properties, we substitute $dN_t^{kj} = dM_t^{kj} + \lambda^{kj}(t, S_t)\, dt$ in (4.11) and rearrange terms to obtain

$$
\theta_t\, dS_t = dv(t, S_t, \eta_t) + J_t\, dt + \sum_{\substack{k,j=1 \\ k \neq j}}^{m} K_t^{kj}\, dM_t^{kj} \quad (4.12)
$$

with bounded integrands

$$
\begin{aligned}
J_t := {} & \delta(t, S_t, \eta_{t-}, v(t, S_t, \eta_{t-})) \\
& + \sum_{\substack{k,j=1 \\ k \neq j}}^{m} \lambda^{kj}(t, S_t) f^{kj}(t, S_t, v(t, S_t, \eta_{t-})) I_{\{k\}}(\eta_{t-}),
\end{aligned}
$$

$$
K_t^{kj} := -(v(t, S_t, j) - v(t, S_t, k)) I_{\{k\}}(\eta_{t-}) I_{\{k \neq j\}}.
$$

For every $p \in [1, \infty)$, $M^{kj}$ is an $\mathcal{H}^p$-martingale under $P$ and $\widehat{P}$. The same property holds for $L$ and for the last term in (4.12) since $v$, $f^{kj}(\cdot, \cdot, v(\cdot, \cdot, \cdot))$ and $K^{kj}$ are all bounded. As $v$ and $J$ are bounded, the right-hand side of (4.12) is the sum of a bounded process and an $\mathcal{H}^p$-martingale under $\widehat{P}$. Hence $\int \theta\, dS$ is well defined on $[0, T]$ and also in $\mathcal{H}^p(\widehat{P})$ since it is already a local $\widehat{P}$-martingale. Because $dP/d\widehat{P}$ has all moments under $\widehat{P}$ due to (4.1), the asserted integrability of $\langle \int \theta\, dS \rangle_T$ under $P$ follows via the Burkholder–Davis–Gundy and Hölder inequalities.

(iii) Denote by $U := \sum_{k \neq j} \int K^{kj}\, dM^{kj}$ the last term from (4.12) so that the local $\widehat{P}$-martingale $\int \theta\, dS$ differs from $U$ only by a bounded process. Both $U$ and $L$ are finite sums of stochastic integrals with respect to $M^{kj}$ of bounded integrands, and so it is enough to prove that each ($P$- and $\widehat{P}$-martingale) $M^{kj}$ is in BMO. But this is clear since $M^{kj}$ has bounded jumps and $\langle M^{kj} \rangle = \int \lambda^{kj}(t, S_t)\, dt$ is bounded for both $P$ and $\widehat{P}$. $\square$

As already mentioned, Theorem 4.4 immediately gives:

COROLLARY 4.5. *Under the assumptions of Theorem 4.4, the strategy $(\theta^0, \theta)$ holding risky assets according to $\theta = (\theta^1, \ldots, \theta^d)$ from (4.10) and investing*

$$\theta_t^0 := v(0, S_0, \eta_0) + L_t + \int_0^t \theta_u \, dS_u - \theta_t S_t, \qquad t \in [0, T],$$

*into the numeraire used for discounting is locally risk-minimizing for $H$. If $\Gamma = 0$, this strategy is even risk-minimizing under $P$.*

PROOF. The first assertion follows from Theorem 3.5 in Schweizer (2001). The second is then clear because local risk-minimization coincides with risk-minimization if $S$ is a local $P$-martingale, that is, for $\Gamma \equiv 0$. □

To obtain a nice representation for the valuation function $v$, we introduce the process

$$Y_t := \widehat{E}[H|\mathcal{F}_t], \qquad t \in [0, T].$$

The structure (4.3) of the claim implies that $H$ has all exponential moments under $\widehat{P}$; this follows from the fact that $h$, $\delta$, $f^{kj}$ are all bounded and $N^{kj}$ has bounded jumps and bounded intensity [see Becherer (2001), Lemma 3.4.1, for details]. Hence $H$ is in $L^p(\widehat{P})$ and $Y$ is an $\mathcal{H}^p$-martingale under $\widehat{P}$ for any $p \in [1, \infty)$.

COROLLARY 4.6. *Under the assumptions of Theorem 4.4, we have the representation*

$$\begin{aligned}
v(t, S_t, \eta_t) &= Y_t - \int_0^t \delta(u, S_u, \eta_{u-}, v(u, S_u, \eta_{u-})) \, du \\
(4.13) &\quad - \int_0^t \sum_{\substack{k,j=1 \\ k \neq j}}^m f^{kj}(u, S_u, v(u, S_u, \eta_{u-})) I_{\{k\}}(\eta_{u-}) \, dN_u^{kj} \\
&= \widehat{E}\bigg[h(S_T, \eta_T) + \int_t^T \delta(u, S_u, \eta_{u-}, v(u, S_u, \eta_{u-})) \, du \\
&\quad + \int_t^T \sum_{\substack{k,j=1 \\ k \neq j}}^m f^{kj}(u, S_u, v(u, S_u, \eta_{u-})) I_{\{k\}}(\eta_{u-}) \, dN_u^{kj} \bigg| \mathcal{F}_t \bigg], \\
&\qquad\qquad\qquad\qquad\qquad\qquad\qquad\qquad\qquad\qquad\qquad t \in [0, T].
\end{aligned}$$

*In particular, the process $V_t := v(t, S_t, \eta_t)$ is the unique solution to the recursive equation (4.7) for $Q = \widehat{P}$, and the function $v$ solves (4.6).*

PROOF. By Theorem 4.4, $v(0, S_0, \eta_0) + \int \theta \, dS + L$ is a $\widehat{P}$-martingale with final value $H$ by (4.5), and so this martingale coincides with $Y$. Hence the first equality in (4.13) follows from (4.11) by integrating, and the second follows from the definitions of $Y$ and $H$. $\square$

The representation of $v$ in (4.13) is given by the conditional expectation under $\widehat{P}$ of future payments. This allows to view $\widehat{P}$ as the pricing measure associated to our approach and thereby confirms (4.6). In the special case when the payoff functions $\delta$ and $f^{kj}$ do not depend on the $v$-argument, (4.13) can also be used to compute $v$ by Monte Carlo methods and thus provides an alternative to the numerical solution of the PDE system (4.8) and (4.9).

REMARK 4.7. To motivate how the PDE (4.8) arises, let us recall from Föllmer and Schweizer (1991) that (4.5) is basically the Galtchouk–Kunita–Watanabe decomposition of $H$ under $\widehat{P}$. For the martingale $Y = \widehat{E}[H|\mathbb{F}]$, the structure of $H$ in (4.3) yields

$$
\begin{aligned}
Y_t = {}& \int_0^t \delta(u, S_u, \eta_{u-}, v(u, S_u, \eta_{u-})) \, du \\
& + \int_0^t \sum_{k \neq j} f^{kj}(u, S_u, v(u, S_u, \eta_{u-})) I_{\{k\}}(\eta_{u-}) \, dN_u^{kj} \\
& + \widehat{E}\bigg[ h(S_T, \eta_T) + \int_t^T \delta(u, S_u, \eta_{u-}, v(u, S_u, \eta_{u-})) \, du \\
& \qquad + \int_t^T \sum_{k \neq j} f^{kj}(u, S_u, v(u, S_u, \eta_{u-})) I_{\{k\}}(\eta_{u-}) \, dN_u^{kj} \bigg| \mathcal{F}_t \bigg].
\end{aligned}
\tag{4.14}
$$

Since $(S, \eta)$ is a Markov process under $\widehat{P}$, the last term is a function of $(t, S_t, \eta_t)$ only, and we call it $v(t, S_t, \eta_t)$. Assuming that $v$ is sufficiently smooth, we can apply Itô's formula, and since $Y$ is a $\widehat{P}$-martingale, all drift terms on the right-hand side must add to 0. Writing out the zero drift condition then produces the PDE system (4.8). However, this is only a heuristic argument because we have given no reason why $v$ should be sufficiently smooth, nor why the last term in (4.14) should be equal to $v(t, S_t, \eta_t)$. The latter is a genuine problem since $v$ appears at many places in (4.14), and this explains why a proof requires a fixed-point argument.

**5. Examples, applications, extensions.** This section outlines some possible applications of our general model. We briefly sketch some links to insurance risk problems and dwell in more detail on issues related to credit risk. There we first explain how the contributions made here fit into an overall perspective, and then show how one can extend the results to incorporate further important aspects in the context of credit risk: hedging by using default-related assets, and stock prices that crash at default.

5.1. *A hybrid reduced form model for credit risk.* In the context of credit risk modeling, it is natural to have a finite-state process $\eta$ representing the rating (or default) state of $\ell$ entities under consideration; we think of $\eta_t \in \{\text{AAA}, \ldots, \text{C}, \text{D (default)}\}^\ell$. As usual in reduced form models, the time of default is defined by a jump of a point process which drives $\eta$. The intensities for rating changes may depend on the stochastic evolution of the Itô process $S$. If the SDE coefficients for $S$ do not depend on $\eta$, then $S$ is a diffusion process, $N$ is a Cox process and $\eta$ is a conditional Markov chain. Then the model falls into the Cox process framework of Lando (1998) which in fact inspired our model. Alternatively, it can be viewed as a (Markovian) reduced form model with state variable process $S$; see Bielecki and Rutkowski (2002) for more explanations and many more references.

REMARK 5.1. One can ask whether tradable asset prices should be affected by the rating from agencies, or rather by some nonobservable solvency state of the firms. We do not address this question here, but use the rating process only as an example for an observable finite-state process which reflects credit risk. One reason to consider the official rating may be the existence of financial products whose payoffs are linked to a rating by agencies.

In comparison to existing credit risk literature, this article makes several contributions. Because of the mutual dependence between assets $S$ and rating process $\eta$, the processes $N$ and $\eta$, conditional on the evolution of $S$, are in general not a Poisson process respectively a Markov chain with intensity matrix $(\lambda^{kj}(t, S_t))_{kj}$. In this regard, the modeling goes beyond the Cox process setting of Lando (1998). Moreover, the model is hybrid in the sense that it incorporates prices $S$ of liquidly traded assets to which the default risk is related, with default intensities being a function of $(t, S_t)$. Such ideas go back to Madan and Unal (1998). If the default intensities of certain firms depend, for instance, on the stock indices of related industries, these indices could be used for hedging and must be taken into account for pricing. Furthermore, we let the payoffs of contingent claims depend on tradable assets and on the default process $\eta$, and we can to some extent also deal with products like convertibles which are related to the stock price of a firm itself; see Section 5.1.2. Yet another contribution is that we consider not only a pure pricing approach under an a priori pricing measure as is usually done in the credit risk literature. We also offer an alternative combined hedging and valuation approach with the aim of minimizing hedging costs. Note that the model should then be specified under the probability measure under which hedging costs are to be minimized. Finally, the solutions to the pricing and hedging problems are described in terms of PDEs under precise conditions; we do not just assume "sufficient regularity" as is often done in the literature.

To exemplify some features of our model, we first consider a single defaultable entity with

$$\eta_t \in \{\text{no default } (\mathfrak{n}), \text{default } (\mathfrak{d})\}$$

and default being an absorbing state. If we choose as discounting numeraire the zero coupon bond with maturity $T$, prices are $T$-forward prices and the payoff of a defaultable zero coupon bond is $h(S_T, \eta_T) := I_{\{\mathfrak{n}\}}(\eta_T)$. Fractional recovery in the case of default by some fraction $R \in [0, 1)$ from treasury or from pre-default (market) value is modeled by choosing $f^{\mathfrak{n}\mathfrak{d}}(t, x, v) := R$ or $f^{\mathfrak{n}\mathfrak{d}}(t, x, v) := Rv$ in the payoff specifications (4.3).

For basket products or counter-party credit risk, one has to deal with multiple default risks and modeling the dependences becomes a crucial task; see Chapter 10 in Schönbucher (2003). In our setting, a natural choice for $\eta$ is the joint rating state of the $\ell$ firms under consideration, that is,

$$(\eta_t^i)_{i=1,\ldots,\ell} = \eta_t \in \{\text{AAA, AA, } \ldots, \text{C, D}\}^\ell, \qquad \ell \in \mathbb{N}.$$

The intensity for a change from the current joint rating state $\eta_{t-}$ to another state $j$ is then

$$\lambda^{kj}(t,x)|_{k=\eta_{t-}, x=S_t} \qquad \text{for } j \in \{\text{AAA, } \ldots, \text{C, D}\}^\ell$$

and can be specified as a function of time $t$, current asset prices $S_t$ and the current overall rating $\eta_{t-}$ itself. This modeling gives a flexible Markovian framework which includes several more specific parametrizations suggested so far. It permits simultaneous defaults of different firms, default correlation because of joint dependence of individual defaults on the common (factor) process $S$, as well as sudden jumps in the default intensities of individual firms when other related firms default. For an implementation, one has to decide on a reduced version with as many or few parameters as can be reasonably fitted from available data.

For a more concrete example, consider $\ell$ default risks and for $\eta$ the simplified state space $\{\mathfrak{n}, \mathfrak{d}\}^\ell$ with no further ratings. We exclude simultaneous defaults and make default an absorbing state for each risk by setting $\lambda^{kj}(t,x) := 0$ when $\sum_{i=1}^\ell I_{\{k^i \neq j^i\}} > 1$ or when $(k^i, j^i) = (\mathfrak{d}, \mathfrak{n})$ for some $i$. Then we are left to choose $\lambda^{kj}(t,x)$ for those $k \neq j$ with $(k^i, j^i) = (\mathfrak{n}, \mathfrak{d})$ for exactly one $i = i^*$ and $k^i = j^i$ elsewhere. As pointed out in Davis and Lo (2001) and Jarrow and Yu (2001), defaults of some firms may trigger jumps in the default intensities of other firms, and a simple way to capture this phenomenon is to take

$$\lambda^{kj}(t,x) = \bar{\lambda}\, a^{(\sum_{i=1}^\ell I_{\{\mathfrak{d}\}}(k^i))} \qquad \text{with } \bar{\lambda} > 0 \text{ and } a \geq 1$$

for those $k \neq j$ with $(k^i, j^i) = (\mathfrak{n}, \mathfrak{d})$ for exactly one $i = i^*$ and $k^i = j^i$ for all other $i$. The basic default intensity of any firm is then $\bar{\lambda}$ and increases by a factor $a \geq 1$ whenever another firm has defaulted. Modeling default dependences between $\ell$ firms is thus reduced to two parameters. To obtain higher default intensities when stock indices are low, practitioners have suggested to model the default intensity as a function of related tradable assets; see Arvanitis and Gregory (2001) or Davis and Lischka (2002). Pushing the example a bit further, one could incorporate such effects in our model by taking $\bar{\lambda} = \bar{\lambda}(t, S_t)$ as a function of a single $(d = 1)$ tradable index $S$. A next step could be a factor model with different groups of firms, and so on.

Instead of considering now one particular application and writing down the resulting PDE system for that problem, we explain in the following two subsections how to extend our framework to address additional aspects of pricing that seem relevant in the context of credit risk.

5.1.1. *Perfect hedging in a completed market.* Up to now, we have assumed that $S$ represents the price processes of all tradable assets which are available for dynamic hedging. But for some applications, it is natural to suppose that there are additional defaultable securities which are liquidly traded, and one feels it should be possible to replicate a given defaultable contingent claim if there are enough tradable securities which are sensitive to the default event. We now show how to make this intuition precise.

To explain the idea in the simplest setting, we consider a single firm which can only be in default or not in default. So we take the model from Section 3.1 for the process $(S, \eta)$ with $S_t \in D = (0, \infty)$ and $\eta_t \in \{\text{no default } (\mathfrak{n}), \text{default } (\mathfrak{d})\}$. The default state is absorbing so that $\lambda^{\mathfrak{dn}}(t, x) := 0$. If $\eta_0 = \mathfrak{n}$ so that the firm is not in default at time $t = 0$, the default time is

$$\tau := \inf\{t \in [0, T] | \eta_t = \mathfrak{d}\} = \inf\{t \in [0, T] | \Delta N_t^{\mathfrak{nd}} \neq 0\}.$$

We suppose that the drift of $S$ vanishes, that is, $\gamma(t, x, k) := 0$ for all $t, x, k$, and we think of $P$ as an a priori pricing measure. But in contrast to Section 4.1, we assume here that there is in addition to $S$ a further tradable security which is sensitive to the default event. For concreteness, we choose as discounting numeraire the (nondefaultable) zero coupon bond with maturity $T$ and suppose that the defaultable zero coupon bond with maturity $T$ and recovery from treasury is tradable. By Theorem 4.4 and Corollary 4.6, the dynamics of its market price process $\bar{Y}_t := E[I_{\{\tau > T\}} + R I_{\{\tau \leq T\}} | \mathcal{F}_t]$ with $R \in [0, 1)$ are given by

$$d\bar{Y}_t = d\bar{v}(t, S_t, \eta_t) = \bar{v}_x(t, S_t, \eta_{t-}) \, dS_t + (\bar{v}(t, S_t, \mathfrak{d}) - \bar{v}(t, S_t, \mathfrak{n})) \, dM_t^{\mathfrak{nd}}$$

and are described by the solution $\bar{v}^k(t, x) = \bar{v}(t, x, k)$ to the PDE system

$$0 = \bar{v}_t^k(t, x) + \tfrac{1}{2} x^2 \sigma^2(t, x, k) \bar{v}_{xx}^k(t, x) + I_{\{k=\mathfrak{n}\}} \lambda^{\mathfrak{nd}}(t, x)(\bar{v}^{\mathfrak{d}}(t, x) - \bar{v}^{\mathfrak{n}}(t, x))$$

with $(t,x) \in [0,T) \times (0,\infty)$, $k = \mathfrak{n}, \mathfrak{d}$ and boundary conditions $\bar{v}^k(T,x) = I_{\{\mathfrak{n}\}}(k) + RI_{\{\mathfrak{d}\}}(k)$. In the present setting, $\bar{v}^\mathfrak{d}(t,x) = R$, and the payoff of a defaultable claim of the form (4.3) is

$$H = h(S_T, \eta_T) + \int_0^T \delta(t, S_t, \eta_{t-}, v(t, S_t, \eta_{t-}))\, dt + I_{\{\tau \leq T\}} f^{\mathfrak{n}\mathfrak{d}}(\tau, S_\tau, v(\tau, S_\tau, \mathfrak{n})).$$

Again by Theorem 4.4 and Corollary 4.6, the dynamics of $Y_t := E[H|\mathcal{F}_t]$, $t \in [0,T]$, are

$$dY_t = v_x(t, S_t, \eta_{t-})\, dS_t + (v(t, S_t, \mathfrak{d}) - v(t, S_t, \mathfrak{n}) + f^{\mathfrak{n}\mathfrak{d}}(t, S_t, v(t, S_t, \mathfrak{n})))\, dM_t^{\mathfrak{n}\mathfrak{d}},$$

where $v^k(t,x) = v(t,x,k)$ denotes the solution to the PDE system

$$0 = v_t^k(t,x) + \tfrac{1}{2}x^2\sigma^2(t,x,k)v_{xx}^k(t,x)$$
$$+ \delta(t,x,k) + I_{\{\mathfrak{n}\}}(k)\lambda^{\mathfrak{n}\mathfrak{d}}(t,x)(v^\mathfrak{d}(t,x) - v^\mathfrak{n}(t,x) + f^{\mathfrak{n}\mathfrak{d}}(t,x,v^\mathfrak{n}(t,x)))$$

with $(t,x) \in [0,T) \times (0,\infty)$, $k = \mathfrak{n}, \mathfrak{d}$ and boundary conditions $v^k(T,x) = h(x,k)$. From the dynamics of $\bar{Y}$ and $Y$, we conclude that

$$dY_t = \left(v_x(t, S_t, \eta_{t-})\right.$$
$$\left. - \frac{v(t, S_t, \mathfrak{d}) - v(t, S_t, \mathfrak{n}) + f^{\mathfrak{n}\mathfrak{d}}(t, S_t, v(t, S_t, \mathfrak{n}))}{\bar{v}(t, S_t, \mathfrak{d}) - \bar{v}(t, S_t, \mathfrak{n})} \bar{v}_x(t, S_t, \eta_{t-})\right) dS_t$$
$$+ \left(\frac{v(t, S_t, \mathfrak{d}) - v(t, S_t, \mathfrak{n}) + f^{\mathfrak{n}\mathfrak{d}}(t, S_t, v(t, S_t, \mathfrak{n}))}{\bar{v}(t, S_t, \mathfrak{d}) - \bar{v}(t, S_t, \mathfrak{n})}\right) d\bar{Y}_t.$$

Hence, the claim $H = Y_T$ can be replicated by self-financing dynamic trading in the risky securities $S$ and $\bar{Y}$, and the hedging strategy is described in terms of the solutions $v$ and $\bar{v}$ to some PDEs. By absence of arbitrage, the value process of $H$ is then given by the wealth process $Y$ of the replicating strategy. This can be viewed as a concrete example (describing the hedging strategy by a PDE) for the more abstract and general representation theorems by Bélanger, Shreve and Wong (2004) and Blanchet-Scalliet and Jeanblanc (2004). For other results on the hedging of *replicable* defaultable claims, see the recent article by Jeanblanc and Rutkowski (2003).

5.1.2. *Stock prices with downward jumps at default.* So far, we have modeled by $S$ the price processes of tradable assets and since $S$ is continuous, we interpret them as stock indices of industries related to the considered defaultable entities. If we consider the tradable stock of a firm itself, things change because we expect this price to jump downwards when the firm defaults on some of its payment obligations. Since products like convertible bonds involve both the stock price and the default state of a firm, this case is very relevant for applications. This section shows how to extend our framework to deal with such situations and obtain similar PDE results as before. We explain the ideas in the simplest setting with a single firm which can be either in default or not, and we restrict ourselves to a pure pricing approach without considering hedging issues here.

We start with the model of Section 3.1 for $(S, \eta)$ with $\eta_t \in \{\mathfrak{n}, \mathfrak{d}\}$ and $S_t \in D = (0, \infty)$. The basic idea is to view $S$ as the *pre-default value* of the stock. The firm is not in default at time $t = 0$ and default is an absorbing state so that $\eta_0 := \mathfrak{n}$ and $\lambda^{\mathfrak{d}\mathfrak{n}}(t, x) := 0$. Default happens at the stopping time

$$(5.1) \qquad \tau := \inf\{t \in [0, T] | \eta_t = \mathfrak{d}\} = \inf\{t \in [0, T] | \Delta N_t^{\mathfrak{n}\mathfrak{d}} \neq 0\},$$

and the events of default or no default up to time $t$ are given by $\{t \geq \tau\} = \{\eta_t = \mathfrak{d}\}$ and $\{t < \tau\} = \{\eta_t = \mathfrak{n}\}$. For a fractional recovery constant $\mathcal{R} \in [0, 1]$, the *firm's stock price* $\bar{S}$ is

$$(5.2) \quad \bar{S}_t := S_t I_{\{t < \tau\}} + \mathcal{R} S_t I_{\{t \geq \tau\}} = S_t - (1 - \mathcal{R}) S_t I_{\{t \geq \tau\}}, \qquad t \in [0, T].$$

Hence the stock price drops to a fraction $\mathcal{R}$ of its pre-default value when default happens. By (5.2) we have $d\bar{S} = \bar{S}_-(dS - (1 - \mathcal{R}) I_{\{\mathfrak{n}\}}(\eta) \, dN^{\mathfrak{n}\mathfrak{d}})$ and by the SDE (3.4) for $S$ we conclude

$$d\bar{S}_t = \bar{S}_{t-}((\gamma(t, S_t, \eta_{t-}) - (1 - \mathcal{R}) I_{\{\mathfrak{n}\}}(\eta_{t-}) \lambda^{\mathfrak{n}\mathfrak{d}}(t, S_t)) \, dt$$
$$+ \sigma(t, S_t, \eta_{t-}) \, dW_t - (1 - \mathcal{R}) I_{\{\mathfrak{n}\}}(\eta_{t-}) \, dM_t^{\mathfrak{n}\mathfrak{d}})$$

with the martingale $M^{\mathfrak{n}\mathfrak{d}}$ as in (4.2); this uses that $S_t = S_{t-} = \bar{S}_{t-}$ on $\{\eta_{t-} = \mathfrak{n}\}$. Since we want $P$ to be a pricing measure, $\bar{S}$ should be a local $P$-martingale, and if we start from a volatility function $\sigma$ and a default intensity function $\lambda^{\mathfrak{n}\mathfrak{d}}$, this requires to take as drift

$$\gamma(t, x, k) := (1 - \mathcal{R}) I_{\{\mathfrak{n}\}}(k) \lambda^{\mathfrak{n}\mathfrak{d}}(t, x).$$

With this choice, the SDE for $\bar{S}$ can be written without reference to $S$ as

$$d\bar{S}_t = \bar{S}_{t-}(\bar{\sigma}(t, \bar{S}_{t-}, \eta_{t-}) \, dW_t - (1 - \mathcal{R}) I_{\{\mathfrak{n}\}}(\eta_{t-}) \, dM_t^{\mathfrak{n}\mathfrak{d}}),$$

with $\bar{\sigma}(t, x, \mathfrak{n}) := \sigma(t, x, \mathfrak{n})$ and $\bar{\sigma}(t, x, \mathfrak{d}) := \sigma(t, x/\mathcal{R}, \mathfrak{d})$ if $\mathcal{R} > 0$ or $\bar{\sigma}(t, x, \mathfrak{d})$ arbitrary if $\mathcal{R} = 0$. It is clear that $\bar{S}$ is then a local $(P, \mathbb{F})$-martingale. To derive the pricing PDE with $P$ as pricing measure, we start for $\mathcal{R} \in [0, 1]$ by considering claims of the form

$$H = \bar{h}(\bar{S}_T) I_{\{\tau > T\}} + \int_0^\tau \bar{\delta}(t, \bar{S}_t, \bar{v}(\tau, \bar{S}_{\tau-})) \, dt + \bar{f}(\tau, \bar{S}_{\tau-}, \bar{v}(\tau, \bar{S}_{\tau-})) I_{\{\tau \leq T\}}$$

with sufficiently regular payoff functions $\bar{h}$, $\bar{\delta}$ and $\bar{f}$ (precise conditions are given below). As before, $\bar{v}$ denotes a valuation function which gives the value of outstanding payments from the claim, given that default has not yet happened, and which still has to be determined. Using definition (5.1) of $\tau$ and the equality $\bar{S}_{t-} = S_{t-} = S_t$ on $[\![0,\tau]\!]$, one can rewrite the claim as

$$H = h(S_T, \eta_T) + \int_0^T \delta(t, S_t, \eta_{t-}, v(t, S_t, \eta_{t-}))\, dt$$

$$+ \int_0^T I_{\{\mathfrak{n}\}}(\eta_{t-}) f^{\mathfrak{n}\mathfrak{d}}(t, S_t, v(t, S_t, \mathfrak{n}))\, dN_t^{\mathfrak{n}\mathfrak{d}}$$

with $h(x,k) = \bar{h}(x) I_{\{k=\mathfrak{n}\}}$, $\delta(t,x,k,v) = \bar{\delta}(t,x,v) I_{\{k=\mathfrak{n}\}}$, $f^{\mathfrak{n}\mathfrak{d}}(t,x,v) = \bar{f}(t,x,v)$ and $v(t,x,k) = \bar{v}(t,x) I_{\{k=\mathfrak{n}\}}$. This represents $H$ in a form like (4.3) which refers to $S$ and $\eta$ instead of $\bar{S}$ and $\tau$. Now suppose that $h$, $\delta$ and $f^{\mathfrak{n}\mathfrak{d}}$ satisfy the regularity conditions stated after (4.3) and recall from Lemma 4.3 (with $Q = P$) that there is at most one valuation process for a claim with a recursive payoff structure like $H$. One can then show that the price $E[H]$ at time 0 equals $v(0, S_0, \mathfrak{n})$, where $v$ is (see Theorem 2.4) the unique bounded solution to the PDE system

$$\begin{aligned}
0 &= v_t(t, x, \mathfrak{n}) + (1-\mathcal{R})\lambda^{\mathfrak{n}\mathfrak{d}}(t, x) x v_x(t, x, \mathfrak{n}) \\
&\quad + \tfrac{1}{2}\sigma^2(t, x, \mathfrak{n}) x^2 v_{xx}(t, x, \mathfrak{n}) + \delta(t, x, \mathfrak{n}, v(t, x, \mathfrak{n})) \\
&\quad + \lambda^{\mathfrak{n}\mathfrak{d}}(t, x)(v(t, x, \mathfrak{d}) - v(t, x, \mathfrak{n}) + f^{\mathfrak{n}\mathfrak{d}}(t, x, v(t, x, \mathfrak{n}))), \\
0 &= v_t(t, x, \mathfrak{d}) + \tfrac{1}{2}\sigma^2(t, x, \mathfrak{d}) x^2 v_{xx}(t, x, \mathfrak{d}),
\end{aligned}$$
(5.3)

for $(t,x) \in [0,T) \times (0,\infty)$, with boundary conditions $v(T,x,k) = h(x,k)$ for $k = \mathfrak{n}, \mathfrak{d}$. Since $v(t,x,\mathfrak{d}) \equiv 0$ is trivial, we are left with one (simplified) PDE for $v(t,x,\mathfrak{n})$. For fractional recovery of the claim in default, that is, $f^{\mathfrak{n}\mathfrak{d}}(t,x,v) = Rv$ with recovery constant $R \in [0,1]$ (for the claim), the last term in (5.3) simplifies further to $-\lambda^{\mathfrak{n}\mathfrak{d}}(t,x)(1-R)v(t,x,\mathfrak{n})$, exhibiting the well-known structure "default intensity times loss fraction." The proof of the above valuation result goes as for Theorem 4.4; one applies Itô's formula to $v(t, S_t, \eta_t)$ to show that

$$Y_t = v(t, S_t, \eta_t) + \int_0^t \delta(u, S_u, \eta_{u-}, v(t, S_u, \eta_{u-}))\, du$$

$$+ \int_0^t I_{\{\mathfrak{n}\}}(\eta_{u-}) f^{\mathfrak{n}\mathfrak{d}}(u, S_u, v(t, S_u, \mathfrak{n}))\, dN_u^{\mathfrak{n}\mathfrak{d}}$$

is a martingale with final value $Y_T = H$ so that $E[H] = Y_0 = v(0, S_0, \eta_0)$.

In the case when $\mathcal{R}$ is strictly positive, the stock price at default cannot drop to zero. By (5.2), we then have not only $S_t = \bar{S}_{t-}$ on $[\![0,\tau]\!]$, but also $S_t = \bar{S}_{t-}/\mathcal{R}$ on $]\!]\tau, T]\!]$. This permits to consider more general claims of the form

$$H = \bar{h}(\bar{S}_T, \eta_T) + \int_0^T \bar{\delta}(t, \bar{S}_{t-}, \eta_{t-}, \bar{v}(t, \bar{S}_{t-}, \eta_{t-}))\, dt$$

$$+ I_{\{\tau \leq T\}} \bar{f}(\tau, \bar{S}_{\tau-}, \bar{v}(\tau, \bar{S}_{\tau-}, \eta_{\tau-}))$$

$$= h(S_T, \eta_T) + \int_0^T \delta(t, S_t, \eta_{t-}, v(t, S_t, \eta_{t-})) \, dt$$

$$+ \int_0^T I_{\{\mathfrak{n}\}}(\eta_{t-}) f^{\mathfrak{n}\mathfrak{d}}(t, S_t, v(t, S_t, \eta_{t-})) \, dN_t^{\mathfrak{n}\mathfrak{d}},$$

where the functions with argument $S$ or $\bar{S}$ are related by $h(x,k) = I_{\{\mathfrak{n}\}}(k)\bar{h}(x,\mathfrak{n}) + I_{\{\mathfrak{d}\}}(k)\bar{h}(x/\mathcal{R},\mathfrak{d})$, $\delta(t,x,k,v) = I_{\{\mathfrak{n}\}}(k)\bar{\delta}(t,x,\mathfrak{n},v) + I_{\{\mathfrak{d}\}}(k)\bar{\delta}(t,x/\mathcal{R},\mathfrak{d},v)$, $f^{\mathfrak{n}\mathfrak{d}}(t,x,v) = \bar{f}(t,x,v)$ and $v(t,x,k) = \bar{v}(t,x,k)I_{\{\mathfrak{n}\}}(k) + \bar{v}(t,x/\mathcal{R},\mathfrak{d})I_{\{\mathfrak{d}\}}(k)$. If $h$, $\delta$ and $f^{\mathfrak{n}\mathfrak{d}}$ satisfy the assumptions after (4.3), arguments analogous to those above yield the price $E[H] = v(0,S_0,\mathfrak{n})$ for $v$ now solving

$$0 = v_t(t,x,\mathfrak{n}) + (1-\mathcal{R})\lambda^{\mathfrak{n}\mathfrak{d}}(t,x)xv_x(t,x,\mathfrak{n})$$
$$+ \tfrac{1}{2}\sigma^2(t,x,\mathfrak{n})x^2 v_{xx}(t,x,\mathfrak{n}) + \delta(t,x,\mathfrak{n},v(t,x,\mathfrak{n}))$$
$$+ \lambda^{\mathfrak{n}\mathfrak{d}}(t,x)(v(t,x,\mathfrak{d}) - v(t,x,\mathfrak{n}) + f^{\mathfrak{n}\mathfrak{d}}(t,x,v(t,x,\mathfrak{n}))),$$
$$0 = v_t(t,x,\mathfrak{d}) + \tfrac{1}{2}\sigma^2(t,x,\mathfrak{d})x^2 v_{xx}(t,x,\mathfrak{d}) + \delta(t,x,\mathfrak{d},v(t,x,\mathfrak{d})),$$

for $(t,x) \in [0,T) \times (0,\infty)$, with boundary conditions $v(T,x,k) = h(x,k)$ for $k = \mathfrak{n}, \mathfrak{d}$. This system is in general more complicated because $v(t,x,\mathfrak{d})$ will no longer vanish if the claim $H$ comprises payments after default time $\tau$.

5.2. *Insurance risk.* A second major area for examples and applications is the treatment of risks at the interface of finance and insurance. This has attracted much attention recently, and we briefly mention a few links to the present work. We concentrate on valuing insurance-related products and leave aside other important topics such as portfolio optimization for insurers or ruin probabilities in the presence of financial markets.

The basic idea is to study a model containing both financial and actuarial components where $S$ describes as before the tradable assets in the financial market, while $\eta$ is now related to the evolution of some insurance contract. For a concrete example, think of $\eta_t$ as the state of health at time $t$ of a policy holder; this may include death if we consider a life insurance policy. The functions $h$, $\delta$ and $f$ specify the contract's final payment $h$ at expiration, continuous payments (e.g., premiums) at rate $\delta$ and lump sum payments $f$ that occur on passing from one state to another. One example of a recursive payoff here arises if death benefits are a fraction of the reserve; see Ramlau-Hansen (1990). In general, the goal is to determine a value for the contract given by $(h, \delta, f)$ and possibly also to hedge it.

While such products per se are not new, recent developments have started to put emphasis on their valuation by market-based methods that go beyond traditional actuarial approaches. Examples include risk-minimization or indifference pricing [see Møller (2001, 2003)] or computation of market-based reserves [see Steffensen (2000)]. In comparison with existing work, we offer here two contributions. We abandon the assumption of independence between the financial and actuarial risk factors imposed in Møller (2001, 2003), and give a rigorous construction for the model with mutual dependences between $S$ and $\eta$. Such dependences can be important for products in non-life insurance; a catastrophic insurance event might, for instance, affect stock indices of related industries. Moreover, we provide precise conditions and existence results for PDEs and fixed points. This contributes to the study of related actuarial applications in Steffensen (2000) and gives a sound mathematical basis to previous intuitive ideas. It will be interesting to see how the methods developed here can be used further, for instance to prove verification results for other problems in the area of insurance.

## APPENDIX

This appendix contains some proofs that were omitted from the main body of the article.

PROOF OF LEMMA 3.2. This is similar to the first step in the proof of Theorem 4.4. We apply Itô's formula to $v(t, S_t, \eta_t)$, substitute $dN_t^{kj} = dM_t^{kj} + \lambda^{kj}(t, S_t)\,dt$ and use the PDE (3.7) to conclude that the drift term vanishes on $[0, T')$. This yields

$$dv(t, S_t, \eta_t) = ((\operatorname{grad}_x v)^{\operatorname{tr}} \Sigma)(t, S_t, \eta_{t-})\,dW_t$$
$$+ \sum_{\substack{j=1 \\ j \neq k}}^{m} \lambda^{kj}(t, x)(v(t, S_t, j) - v(t, S_t, k)) I_{\{k\}}(\eta_{t-})\,dM_t^{kj},$$
$$t \in [0, T').$$

Since $v$ and $\lambda^{kj}$ are bounded and the martingales $M^{kj}$ are in $\mathcal{H}^p(P)$ for any $p \in [1, \infty)$ by the remark following (4.2), all stochastic integrals with respect to $M^{kj}$ are in $\mathcal{H}^p(P)$, and so this must hold for the one with respect to $W$ as well. In particular, the above equation extends to $[0, T']$ and the process $v(t, S_t, \eta_t)$, $t \in [0, T']$, is a martingale. □

PROOF OF LEMMA 4.3. Suppose $V^1$ and $V^2$ are bounded semimartingales which both satisfy (4.7). By using the martingale compensator $\lambda^{kj}(t, S_t)\,dt$ for the jump process $N^{kj}$ [which is the same under $P$ and $\widehat{P}$, see (4.2)], the recursive representation (4.7) of $V^i$ can be rewritten as

$$(A.1) \qquad V_t^i = E_Q\bigg[h(S_T, \eta_T) + \int_t^T g(u, S_u, \eta_{u-}, V_{u-}^i)\,du\Big|\mathcal{F}_t\bigg],$$
$$t \in [0, T], i = 1, 2,$$

with $g(t,x,k,y) := \delta(t,x,k,y) + \sum_{j\neq k} \lambda^{kj}(t,x) f^{kj}(t,x,y)$. Since all $\delta(\cdot,\cdot,k,\cdot)$ and $f^{kj}$ are locally Lipschitz in $y$, uniformly in $(t,x)$, $g$ is also locally Lipschitz in $y$, uniformly in $(t,x,k)$. If $K < \infty$ is an upper bound for $|V^1|, |V^2|$, the local Lipschitz property gives an $L = L(K) < \infty$ such that $|g(t,x,k,y^1) - g(t,x,k,y^2)| \leq L|y^1 - y^2|$ for any $|y^i| \leq K$, $i=1,2$, and all $(t,x,k)$. Similarly as for Proposition 2.1, it follows from (A.1) that we have for any $\beta \in (0,\infty)$

$$
\begin{aligned}
&e^{-\beta(T-t)}|V_t^1 - V_t^2| \\
&\leq e^{-\beta(T-t)} E_Q\left[\int_t^T |g(u, S_u, \eta_{u-}, V_{u-}^1) - g(u, S_u, \eta_{u-}, V_{u-}^2)|\, du \Big| \mathcal{F}_t\right] \\
&\leq e^{-\beta(T-t)} \left\|\sup_{u\in[0,T]} e^{-\beta(T-u)} L |V_u^1 - V_u^2|\right\|_{L^\infty} \int_t^T e^{+\beta(T-u)}\, du \\
&\leq \frac{L}{\beta} \left\|\sup_{u\in[0,T]} e^{-\beta(T-u)}|V_u^1 - V_u^2|\right\|_{L^\infty}, \qquad t \in [0,T].
\end{aligned}
$$

So if the RCLL semimartingales $V^1$ and $V^2$ were not indistinguishable, considering the supremum of the above inequality over $t \in [0,T]$ would lead to a contradiction for $\beta > L$. □

The second part of the above proof varies an argument from Duffie and Epstein (1992) (Appendix with Skiadas). Their result cannot be applied directly because the function $g$ above is only locally Lipschitz. And the exponential weighting permits to shorten the proof.

**Acknowledgments.** The first author thanks Fredi Tröltzsch, TU Berlin, and participants of workshops in Copenhagen and Oberwolfach for discussions. Earlier versions contained an incomplete proof of Proposition 2.3, and we are grateful to Gallus Steiger for pointing out this gap. All remaining errors are ours.

Department of Mathematics
Imperial College London
London SW7 2AZ
United Kingdom
e-mail: dirk.becherer@imperial.ac.uk

Department of Mathematics
ETH Zürich
CH-8092 Zürich
Switzerland
e-mail: mschweiz@math.ethz.ch